\documentclass[11pt,oneside]{amsart}

\usepackage{tikz-cd}

\usepackage{newtxtext}
\usepackage{newtxmath}
\usepackage{fix-cm}

\DeclareMathAlphabet{\mathcal}{OMS}{cmsy}{m}{n}

\usepackage[arrow,curve,matrix,tips,2cell]{xy}

\usepackage[hmargin=2cm,vmargin=1.5cm]{geometry}

\usepackage{bm}

\usepackage{fancyhdr}

\usepackage[obeyspaces]{url}

\usepackage{amsmath}
\usepackage{amscd}
\usepackage{amssymb}
\usepackage{latexsym}
\usepackage{url}

\usepackage{graphicx}

\newtheorem{theorem}{Theorem}[section]
\newtheorem*{theorem*}{Theorem}
\newtheorem{lemma}[theorem]{Lemma}
\newtheorem*{lemma*}{Lemma}
\newtheorem{corollary}[theorem]{Corollary}
\newtheorem{proposition}[theorem]{Proposition}

\newtheorem{remark}[theorem]{Remark}
\newtheorem{definition}[theorem]{Definition}
\newtheorem*{definition*}{Definition}

\newtheorem{question}[theorem]{Question}
\newtheorem*{question*}{Question}

\newtheorem{example}[theorem]{Example}
\newtheorem{examples}[theorem]{Examples}

\makeatletter
\def\revddots{\mathinner{\mkern1mu\raise\p@
\vbox{\kern7\p@\hbox{.}}\mkern2mu
\raise4\p@\hbox{.}\mkern2mu\raise7\p@\hbox{.}\mkern1mu}}
\makeatother 
\newcommand{\bgl}{\begin{equation}} 
\newcommand{\egl}{\end{equation}}
\newcommand{\bgloz}{\begin{equation*}} 
\newcommand{\egloz}{\end{equation*}}
\newcommand{\bgln}{\begin{eqnarray}} 
\newcommand{\egln}{\end{eqnarray}}
\newcommand{\bglnoz}{\begin{eqnarray*}} 
\newcommand{\eglnoz}{\end{eqnarray*}}
\newcommand{\btheo}{\begin{theorem}}
\newcommand{\etheo}{\end{theorem}}
\newcommand{\btheooz}{\begin{theorem*}}
\newcommand{\etheooz}{\end{theorem*}}

\newcommand{\blemma}{\begin{lemma}}
\newcommand{\elemma}{\end{lemma}}
\newcommand{\blemmaoz}{\begin{lemma*}}
\newcommand{\elemmaoz}{\end{lemma*}}
\newcommand{\bproof}{\begin{proof}}
\newcommand{\eproof}{\end{proof}}
\newcommand{\bbew}{\begin{beweis}}
\newcommand{\ebew}{\end{beweis}}
\newcommand{\bremark}{\begin{remark}\em}
\newcommand{\eremark}{\end{remark}}
\newcommand{\bdefin}{\begin{definition}}
\newcommand{\edefin}{\end{definition}}
\newcommand{\bdefinoz}{\begin{definition*}}
\newcommand{\edefinoz}{\end{definition*}}
\newcommand{\bex}{\begin{example}}
\newcommand{\eex}{\end{example}}
\newcommand{\bexs}{\begin{examples}}
\newcommand{\eexs}{\end{examples}}

\newcommand{\bprop}{\begin{proposition}}
\newcommand{\eprop}{\end{proposition}}
\newcommand{\bcor}{\begin{corollary}}
\newcommand{\ecor}{\end{corollary}}
\newcommand{\bfa}{\begin{cases}} 
\newcommand{\efa}{\end{cases}}

\newcommand{\bquestion}{\begin{question}}
\newcommand{\equestion}{\end{question}}
\newcommand{\bquestionoz}{\begin{question*}}
\newcommand{\equestionoz}{\end{question*}}
\newcommand{\cC}{\mathcal C}

\newcommand{\cF}{\mathcal F}

\newcommand{\cK}{\mathcal K}

\newcommand{\cO}{\mathcal O}
\newcommand{\cP}{\mathcal P}

\newcommand{\cR}{\mathcal R}

\newcommand{\cT}{\mathcal T}

\newcommand{\cW}{\mathcal W}
\newcommand{\cX}{\mathcal X}
\newcommand{\cY}{\mathcal Y}
\newcommand{\cZ}{\mathcal Z}

\def\Cz{\mathbb{C}}

\def\Nz{\mathbb{N}}

\def\Tz{\mathbb{T}}
\def\Zz{\mathbb{Z}}

\def\1z{\mathbb{1}}
\newcommand{\an}[1]{``#1''} 
\newcommand{\ti}{\tilde}

\newcommand{\lori}{\longrightarrow}

\newcommand{\ma}{\mapsto} 
\newcommand{\onto}{\twoheadrightarrow} 
\newcommand{\into}{\hookrightarrow} 
\newcommand{\Rarr}{\Rightarrow} 

\def\SEMI{\mbox{$\times\kern-2pt\vrule height5pt width.6pt \kern3pt $}}

\newcommand{\id}{{\rm id}}

\renewcommand{\dim}{{\rm dim}\,}
\newcommand{\rk}{{\rm rk}\,}
\newcommand{\img}{{\rm im\,}}

\newcommand{\ev}{\operatorname{ev}} 
\newcommand{\lspan}{{\rm span}} 
\newcommand{\clspan}{\overline{\lspan}} 
\newcommand{\norm}[1]{\left\|#1\right\|} 
\newcommand{\defeq}{\mathrel{:=}} 
\newcommand{\eqdef}{\mathrel{=:}} 

\newcommand{\dop}{\text{: }} 
\newcommand{\ilim}{\varinjlim} 
\newcommand{\plim}{\varprojlim} 
\newcommand{\supp}{{\rm supp}} 
\newcommand{\lge}{\left\{} 
\newcommand{\rge}{\right\}} 
\newcommand{\lru}{\left(} 
\newcommand{\rru}{\right)} 
\newcommand{\rukl}[1]{\lru #1 \rru} 
\newcommand{\gekl}[1]{\lge #1 \rge} 
\newcommand{\menge}[2]{\gekl{ #1 \dop #2 }} 
\def\bf1{\mathbf{1}}
\newcommand{\dom}{{\rm dom\,}}
\begin{document}

\title{Every classifiable simple C*-algebra has a Cartan subalgebra}

\thispagestyle{fancy}

\begin{abstract}
We construct Cartan subalgebras in all classifiable stably finite C*-algebras. Together with known constructions of Cartan subalgebras in all UCT Kirchberg algebras, this shows that every classifiable simple C*-algebra has a Cartan subalgebra.
\end{abstract}

\author{Xin Li}

\address{Xin Li, School of Mathematical Sciences, Queen Mary University of London, Mile End Road, London E1 4NS, United Kingdom}
\email{xin.li@qmul.ac.uk}


\maketitle


\setlength{\parindent}{0cm} \setlength{\parskip}{0.5cm}

\section{Introduction}
\label{s:Intro}

Classification of C*-algebras has seen tremendous advances recently. In the unital case, the classification of unital separable simple nuclear $\cZ$-stable C*-algebras satisfying the UCT is by now complete. This is the culmination of work by many mathematicians. The reader may consult \cite{KP, Phi, GLN, EGLN, TWW} and the references therein. In the stably projectionless case, classification results are being developed (see \cite{EN,EGLN17a,EGLN17b,GLI,GLII}). It is expected that -- once the stably projectionless case is settled -- the final result will classify all separable simple nuclear $\cZ$-stable C*-algebras satisfying the UCT by their Elliott invariants. This class of C*-algebras is what we refer to as \an{classifiable C*-algebras}. 
\setlength{\parindent}{0.5cm} \setlength{\parskip}{0cm}

To complete these classification results, it is important to construct concrete models realizing all possible Elliott invariants by classifiable C*-algebras. Such models have been constructed -- in the greatest possible generality -- in \cite{Ell} (see also \cite{Tho} which covers special cases). In the stably finite unital case, the reader may also find such range results in \cite{GLN}, where the construction follows the ideas in \cite{Ell} (with slight modifications, so that the models belong to the special class considered in \cite{GLN}). In the stably projectionless case, models have been constructed in a slightly different way in \cite{GLII} (again to belong to the special class of algebras considered) under the additional assumption of a trivial pairing between K-theory and traces.
\setlength{\parindent}{0cm} \setlength{\parskip}{0.5cm}

Recently, the notion of Cartan subalgebras in C*-algebras \cite{Kum,Ren} has attracted attention, due to connections to topological dynamics \cite{Li16,Li17,Li_DQH} and the UCT question \cite{BL16,BL17}. In particular the reformulation of the UCT question in \cite{BL16,BL17} raises the following natural question (see \cite[Question~5.9]{LR}, \cite[Question~16]{Banff17} and \cite[Problems~1 and 2]{OWR17}):

\bquestion
\label{?}
Which classifiable C*-algebras have Cartan subalgebras?
\equestion 
\setlength{\parindent}{0cm} \setlength{\parskip}{0cm}

By \cite{Kum,Ren}, we can equally well ask for groupoid models for classifiable C*-algebras. In the purely infinite case, groupoid models and hence Cartan subalgebras have been constructed in \cite{Spi} (see also \cite[\S~5]{LR}). For special classes of stably finite unital C*-algebras, groupoid models have been constructed in \cite{DPS,Put} using topological dynamical systems. Using a new approach, the goal of this paper is to answer Question~\ref{?} by constructing Cartan subalgebras in all the C*-algebra models constructed in \cite{Ell,GLN,GLII}, covering all classifiable stably finite C*-algebras, in particular in all classifiable unital C*-algebras. Generally speaking, Cartan subalgebras allow us to introduce ideas from geometry and dynamical systems to the study of C*-algebras. More concretely, in view of \cite{BL16,BL17}, we expect that our answer to Question~\ref{?} will lead to progress on the UCT question.
\setlength{\parindent}{0cm} \setlength{\parskip}{0.5cm}

The following two theorems are the main results of this paper. The reader may consult \cite{Kum,Ren} for the definition of twisted groupoids and their relation to Cartan subalgebras, and \cite[\S~2.2]{Ror}, \cite[\S~8.4]{Lin}, \cite{GLN,GLI,GLII} for the precise definition of the Elliott invariant.

\btheo[unital case]
\label{THM1}
Given
\setlength{\parindent}{0cm} \setlength{\parskip}{0cm}

\begin{itemize}
\item a weakly unperforated, simple scaled ordered countable abelian group $(G_0,G_0^+,u)$,
\item a non-empty metrizable Choquet simplex $T$,
\item a surjective continuous affine map $r: \: T \to S(G_0)$,
\item a countable abelian group $G_1$,
\end{itemize}

there exists a twisted groupoid $(G,\Sigma)$ such that
\begin{itemize}
\item $G$ is a principal \'{e}tale second countable locally compact Hausdorff groupoid,
\item $C^*_r(G,\Sigma)$ is a simple unital C*-algebra which can be described as the inductive limit of subhomogeneous C*-algebras whose spectra have dimension at most $3$,
\item the Elliott invariant of $C^*_r(G,\Sigma)$ is given by
$$
  (K_0(C^*_r(G,\Sigma)), K_0(C^*_r(G,\Sigma))^+, [1_{C^*_r(G,\Sigma)}], T(C^*_r(G,\Sigma)), r_{C^*_r(G,\Sigma)}, K_1(C^*_r(G,\Sigma))) \cong (G_0,G_0^+,u,T,r,G_1).
$$
\end{itemize}
\etheo

\btheo[stably projectionless case]
\label{THM2}
Given
\setlength{\parindent}{0cm} \setlength{\parskip}{0cm}

\begin{itemize}
\item countable abelian groups $G_0$ and $G_1$,
\item a non-empty metrizable Choquet simplex $T$,
\item a homomorphism $\rho: \: G_0 \to {\rm Aff}(T)$ which is weakly unperforated in the sense that for all $g \in G_0$, there is $\tau \in T$ with $\rho(g)(\tau) = 0$
\end{itemize}

there exists a twisted groupoid $(G,\Sigma)$ such that
\begin{itemize}
\item $G$ is a principal \'{e}tale second countable locally compact Hausdorff groupoid,
\item $C^*_r(G,\Sigma)$ is a simple stably projectionless C*-algebra with continuous scale in the sense of \cite{Lin91,Lin04,GLI,GLII} which can be described as the inductive limit of subhomogeneous C*-algebras whose spectra have dimension at most $3$,
\item the Elliott invariant of $C^*_r(G,\Sigma)$ is given by
$$
  (K_0(C^*_r(G,\Sigma)), K_0(C^*_r(G,\Sigma))^+, T(C^*_r(G,\Sigma)), \rho_{C^*_r(G,\Sigma)}, K_1(C^*_r(G,\Sigma))) \cong (G_0,\gekl{0},T,\rho,G_1).
$$
\end{itemize}
\etheo
\setlength{\parindent}{0cm} \setlength{\parskip}{0cm}

The condition on $\rho$ means that the pairing between K-theory and traces is weakly unperforated, in the sense of \cite{Ell}. It has been shown in \cite[\S~A.1]{EGLN17b} that this condition of weak unperforation is necessary in the classifiable setting (i.e., it follows from finite nuclear dimension, or $\cZ$-stability).
\setlength{\parindent}{0cm} \setlength{\parskip}{0.5cm}

It is worth pointing out that in the main theorems, the twisted groupoids are constructed explicitly in such a way that the inductive limit structure with subhomogeneous building blocks will become visible at the groupoid level.

\bremark
The original building blocks in \cite{Ell} have spectra with dimension at most two. The reason three-dimensional spectra are needed in this paper is because it is not clear how to realize all possible connecting maps at the level of $K_1$ by Cartan-preserving homomomorphisms using the building blocks in \cite{Ell}. Therefore, the building blocks have to be modified (see \S~\ref{s:AH}). Roughly speaking, the idea is to realize all possible connecting maps in $K_1$ at the level of topological spaces. This however requires three-dimensional spectra because \an{nice} topological spaces (say CW-complexes) of dimension two or lower have torsion-free $K^1$ (because cohomology is torsion-free in all odd degrees for these spaces). The dimension can be reduced to two if $K_1$ is torsion-free (see Corollary~\ref{C:dim} and Remark~\ref{R:S3S1}).
\eremark

In particular, together with the classification results in \cite{KP,Phi,GLN,EGLN,TWW}, the groupoid models in \cite{Spi}, and \cite[Theorem~3.1]{BL16}, we obtain the following
\bcor
A unital separable simple C*-algebra with finite nuclear dimension has a Cartan subalgebra if and only if it satisfies the UCT.
\ecor
\setlength{\parindent}{0cm} \setlength{\parskip}{0cm}
The only reason we restrict to the unital case here is that classification in the stably projectionless case has not been completed yet.
\setlength{\parindent}{0cm} \setlength{\parskip}{0.5cm}

The constructions of the twisted groupoids in Theorems~\ref{THM1} and \ref{THM2} yield the following direct consequences:
\bcor
\label{C:unbdd1}
In the situation of Theorem~\ref{THM1}, suppose that in addition to $(G_0, G_0^+, u)$, $T$, $r$ and $G_1$, we are given a topological cone $\ti{T}$ with base $T$ and a lower semicontinuous affine map $\ti{\gamma}: \: \ti{T} \to [0,\infty]$. Then there exists a twisted groupoid $(\ti{G},\ti{\Sigma})$ such that 

\begin{itemize}
\item $\ti{G}$ is a principal \'{e}tale second countable locally compact Hausdorff groupoid,
\item $C^*_r(\ti{G},\ti{\Sigma})$ is a non-unital hereditary sub-C*-algebra of $C^*_r(G,\Sigma) \otimes \cK$,
\item the Elliott invariant of $C^*_r(\ti{G},\ti{\Sigma})$ is given by
$$
  (K_0(C^*_r(\ti{G},\ti{\Sigma})), K_0(C^*_r(\ti{G},\ti{\Sigma}))^+, \ti{T}(C^*_r(G,\Sigma)), \Sigma_{C^*_r(\ti{G},\ti{\Sigma})}, r_{C^*_r(\ti{G},\ti{\Sigma})}, K_1(C^*_r(\ti{G},\ti{\Sigma}))) \cong (G_0,G_0^+,\ti{T},\ti{\gamma},r,G_1).
$$
\end{itemize}
\ecor

\bcor
\label{C:unbdd2}
In the situation of Theorem~\ref{THM2}, suppose that in addition to $G_0$, $G_1$, $T$ and $\rho$, we are given a topological cone $\ti{T}$ with base $T$ and a lower semicontinuous affine map $\ti{\gamma}: \: \ti{T} \to [0,\infty]$. Then there exists a twisted groupoid $(\ti{G},\ti{\Sigma})$ such that
\setlength{\parindent}{0cm} \setlength{\parskip}{0cm}

\begin{itemize}
\item $\ti{G}$ is a principal \'{e}tale second countable locally compact Hausdorff groupoid,
\item $C^*_r(\ti{G},\ti{\Sigma})$ is a hereditary sub-C*-algebra of $C^*_r(G,\Sigma) \otimes \cK$,
\item the Elliott invariant of $C^*_r(\ti{G},\ti{\Sigma})$ is given by
$$
  (K_0(C^*_r(\ti{G},\ti{\Sigma})), K_0(C^*_r(\ti{G},\ti{\Sigma}))^+, \ti{T}(C^*_r(\ti{G},\ti{\Sigma})), \Sigma_{C^*_r(\ti{G},\ti{\Sigma})}, \rho_{C^*_r(\ti{G},\ti{\Sigma})}, K_1(C^*_r(\ti{G},\ti{\Sigma})) \cong (G_0,\gekl{0},\ti{T},\ti{\gamma},\rho,G_1).
$$
\end{itemize}
\ecor
Note that all the groupoids in Theorems~\ref{THM1}, \ref{THM2} and Corollaries~\ref{C:unbdd1}, \ref{C:unbdd2} are necessarily minimal and amenable. Theorem~\ref{THM1} and Corollary~\ref{C:unbdd1}, together with \cite{Spi}, imply that every classifiable C*-algebra which is not stably projectionless has a Cartan subalgebra. Once the classification of stably projectionless C*-algebras is completed, Theorem~\ref{THM2} and Corollary~\ref{C:unbdd2} will imply that every classifiable stably projectionless C*-algebra has a Cartan subalgebra. Actually, using $\cZ$-stability, we see that all of the above-mentioned classifiable C*-algebras have infinitely many non-isomorphic Cartan subalgebras (compare \cite[Proposition~5.1]{LR}). Moreover, the constructions in this paper show that in every classifiable stably finite C*-algebra, we can even find C*-diagonals (and even infinitely many non-isomorphic ones).
\setlength{\parindent}{0cm} \setlength{\parskip}{0.5cm}

Moreover, more can be said about the twist, and also about the dimension of the spectra of our Cartan subalgebras.
\bcor
\label{C:dim}
The twisted groupoids $(G,\Sigma)$ constructed in the proofs of Theorems~\ref{THM1} and \ref{THM2} have the following additional properties:
\setlength{\parindent}{0cm} \setlength{\parskip}{0cm}

\begin{enumerate}
\item[(i)] If $G_0$ is torsion-free, then the twist $\Sigma$ is trivial, i.e., $\Sigma = \Tz \times G$.
\item[(ii)] If $G_1$ has torsion, then $C^*_r(G,\Sigma)$ is an inductive limit of subhomogeneous C*-algebras whose spectra are three-dimensional, and $\dim(G^{(0)}) = 3$.
\item[(iii)] If $G_1$ is torsion-free and $G_0$ has torsion, $C^*_r(G,\Sigma)$ is an inductive limit of subhomogeneous C*-algebras whose spectra are two-dimensional, and $\dim(G^{(0)}) = 2$.
\item[(iv)] If both $G_0$ and $G_1$ are torsion-free with $G_1 \ncong \gekl{0}$, then $C^*_r(G,\Sigma)$ is an inductive limit of subhomogeneous C*-algebras whose spectra are one-dimensional, and $\dim(G^{(0)}) = 1$.
\item[(v)] If $G_0$ is torsion-free and $G_1 \cong \gekl{0}$, then $C^*_r(G,\Sigma)$ is an inductive limit of one-dimensional non-commutative finite CW-complexes, with $\dim(G^{(0)}) \leq 1$ in Theorem~\ref{THM1} and $\dim(G^{(0)}) = 1$ in Theorem~\ref{THM2}.
\end{enumerate}
\ecor
\setlength{\parindent}{0cm} \setlength{\parskip}{0cm}

In particular, Corollary~\ref{C:dim} implies the following:
\bcor
\label{C:ZWZ}
The Jiang-Su algebra $\cZ$, the Razak-Jacelon algebra $\cW$ and the stably projectionless version $\cZ_0$ of the Jiang-Su algebra of \cite[Definition~7.1]{GLII} have C*-diagonals with one-dimensional spectra. The corresponding twisted groupoids $(G,\Sigma)$ can be chosen so that $\Sigma$ is trivial, i.e., $\Sigma = \Tz \times G$.
\ecor
Concrete groupoid models for $\cZ$, $\cW$ and $\cZ_0$ are described in \S~\ref{s:ZWZ}. It is worth pointing out that a groupoid model has been constructed for $\cZ$ in \cite{DPS} using a different construction (but the precise dimension of the unit space has not been determined in \cite{DPS}). Moreover, G. Szab\'{o} and S. Vaes independently found groupoid models for $\cW$, again using constructions different from ours. Furthermore, independently from \cite{BL17} and the present paper, similar tools to the ones in \cite[\S~3]{BL17} were developed in \cite{AM}, which give rise to groupoid models for $\cZ$ and $\cW$ as well as other examples.
\setlength{\parindent}{0cm} \setlength{\parskip}{0.5cm}

The key tool for all the results in this paper is an improved version of \cite[Theorem~3.6]{BL17}, which allows us to construct Cartan subalgebras in inductive limit C*-algebras. The C*-algebraic formulation reads as follows.
\btheo
\label{THM:CartanLim}
Let $(A_n,B_n)$ be Cartan pairs with normalizers $N_n \defeq N_{A_n}(B_n)$ and faithful conditional expectations $P_n: \: A_n \onto B_n$. Let $\varphi_n: \: A_n \to A_{n+1}$ be injective *-homomorphisms with $\varphi_n(B_n) \subseteq B_{n+1}$, $\varphi_n(N_n) \subseteq N_{n+1}$ and $P_{n+1} \circ \varphi_n = \varphi_n \circ P_n$ for all $n$. Then $\ilim \gekl{B_n;\varphi_n}$ is a Cartan subalgebra of $\ilim \gekl{A_n;\varphi_n}$.
\setlength{\parindent}{0.5cm} \setlength{\parskip}{0cm}

If all $B_n$ are C*-diagonals, then $\ilim \gekl{B_n;\varphi_n}$ is a C*-diagonal of $\ilim \gekl{A_n;\varphi_n}$.
\etheo
\setlength{\parindent}{0cm} \setlength{\parskip}{0cm}
A special case of this theorem is proved in \cite{DP}.

\setlength{\parindent}{0cm} \setlength{\parskip}{0.5cm}

Actually, in addition to Theorem~\ref{THM:CartanLim}, much more is accomplished: Groupoid models are developed for *-homomorphisms such as $\varphi_n$, and the twisted groupoid corresponding to $\rukl{\ilim \gekl{A_n;\varphi_n}, \ilim \gekl{B_n;\varphi_n}}$ as in Theorem~\ref{THM:CartanLim} is described explicitly. These results (in \S~\ref{s:LimCartan}) might be of independent interest.
\setlength{\parindent}{0cm} \setlength{\parskip}{0.5cm}

Applications of these explicit descriptions of groupoid models (for homomorphisms and Cartan pairs) and Theorem~\ref{THM:CartanLim} include a unified approach to Theorems~\ref{THM1}, \ref{THM2}, and explicit constructions of the desired twisted groupoids. The strategy is as follows: C*-algebras with prescribed Elliott invariant have been constructed in \cite{Ell} (see also \cite[\S~13]{GLN} for the unital case). These C*-algebras have all the desired properties as in Theorems~\ref{THM1} and \ref{THM2} and are constructed as inductive limits of subhomogeneous C*-algebras. However, the connecting maps in \cite{Ell} and \cite[\S~13]{GLN} do not preserve the canonical Cartan subalgebras in these building blocks in general. Therefore, a careful choice or modification of the building blocks and connecting maps in the constructions in \cite{Ell,GLN} is necessary in order to allow for an application of Theorem~\ref{THM:CartanLim}. The modification explained in Remark~\ref{R:unitaries} is particularly important. Actually, a more general result is established in \S~\ref{ss:SummConst}, where a class of inductive limits of subhomogeneous C*-algebras is identified, which encompasses all the C*-algebras in Theorems~\ref{THM1}, \ref{THM2} and Corollaries~\ref{C:unbdd1}, \ref{C:unbdd2}, where we can apply Theorem~\ref{THM:CartanLim}.

I am grateful to the organizers Sel{\c{c}}uk Barlak, Wojciech Szyma\'{n}ski and Wilhelm Winter of the Oberwolfach Mini-Workshop \an{MASAs and Automorphisms of C*-Algebras} for inviting me, and for the discussions in Oberwolfach with Sel{\c{c}}uk Barlak which eventually led to this paper. I also thank Sel{\c{c}}uk Barlak and G\'{a}bor Szab\'{o} for helpful comments on earlier drafts. Moreover, I would like to thank the referee for very helpful comments which led to an improved version of Theorem~\ref{THM2} (previous versions of this theorem only covered classifiable stably projectionless C*-algebras with trivial pairing between K-theory and traces).

\section{The constructions of Elliott and Gong-Lin-Niu}
\setlength{\parindent}{0cm} \setlength{\parskip}{0cm}

Let us briefly recall the constructions in \cite{Ell} (see also \cite{EV} for simplifications and further explanations) and \cite[\S~13]{GLN}.

\subsection{The unital case}
\label{ss:UnitalCase}

Let us describe the construction in \cite[\S~13]{GLN}, which is based on \cite{Ell} (with slight modifications). Given $(G_0,G_0^+,u,T,r,G_1)$ as in Theorem~\ref{THM1}, write $G = G_0$, $K = G_1$, and let $\rho: \: G \to {\rm Aff}(T)$ be the dual map of $r$. Choose a dense subgroup $G' \subseteq {\rm Aff}(T)$. Set $H \defeq G \oplus G'$, 
$$
  H^+ \defeq \gekl{(0,0)} \cup \menge{(g,f) \in G \oplus G'}{\rho(g)(\tau) + f(\tau) > 0 \ \forall \ \tau \in T},
$$
and view $u$ in $G$ as an element of $H$. Then $(H,H^+,u)$ becomes a simple ordered group, inducing the structure of a dimension group on $H / {\rm Tor}(H)$. Now construct a commutative diagram
$$
  \xymatrix{
  G_1 \ar[d] \ar[r] & \dotso \ar[r] & G_n \ar[d] \ar[r] & G_{n+1} \ar[d] \ar[r] & & \dotso & \ar[r] & G \ar[d]
  \\
  H_1 \ar[d] \ar[r]^{\gamma_1} & \dotso \ar[r]^{\gamma_{n-1}} & H_n \ar[d] \ar[r]^{\gamma_n} & H_{n+1} \ar[d] \ar[r]^{\gamma_{n+1}} & & \dotso & \ar[r] & H \ar[d]
  \\
  H_1 / G_1 \ar[r] & \dotso \ar[r] & H_n / G_n \ar[r] & H_{n+1} / G_{n+1} \ar[r] & & \dotso & \ar[r] & H / G
  }
$$
where:
\setlength{\parindent}{0cm} \setlength{\parskip}{0cm}

\begin{itemize}
\item $H_n$ is a finitely generated abelian group with $H_n = \bigoplus_i H_n^i$, where for one distinguished index $\bm{i}$, $H_n^{\bm{i}} = \Zz \oplus {\rm Tor}(H_n)$, and for all other indices, $H_n^i = \Zz$;
\item with $(H_n^{\bm{i}})^+ \defeq \gekl{(0,0)} \cup (\Zz_{>0} \oplus {\rm Tor}(H_n))$, $(H_n^i)^+ \defeq \Zz_{\geq 0}$ for all $i \neq \bm{i}$, $H_n^+ \defeq \bigoplus_i (H_n^i)^+ \subseteq H_n^{\bm{i}} \oplus \bigoplus_{i \neq \bm{i}} H_n^i = H_n$ and $u_n = (([n,\bm{i}],\tau_n),([n,i])_{i \neq \bm{i}}) \in H_n^+$, we have 
\begin{equation}
\label{e:H=H}
\ilim \gekl{(H_n,H_n^+,u_n); \gamma_n} \cong (H,H^+,u);
\end{equation}
\item with $G_n \defeq (\gamma_n^{\infty})^{-1}(G)$, where $\gamma_n^{\infty}: \: H_n \to H$ is the map provided by \eqref{e:H=H}, and $G_n^+ \defeq G_n \cap H_n^+$, we have $u_n \in G_n \subseteq H_n$, and \eqref{e:H=H} induces $\ilim \gekl{(G_n,G_n^+,u_n); \gamma_n} \cong (G,G^+,u)$;
\item the vertical maps are the canonical ones.
\end{itemize}
\setlength{\parindent}{0cm} \setlength{\parskip}{0.5cm}

Let $\hat{\gamma}_n: \: H_n / {\rm Tor}(H_n) \eqdef \hat{H}_n = \bigoplus_i \hat{H}_n^i \to \bigoplus_j \hat{H}_{n+1}^j = \hat{H}_{n+1} \defeq H_{n+1} / {\rm Tor}(H_{n+1})$ be the homomorphism induced by $\gamma_n$, where $\hat{H}_n^i = \Zz = \hat{H}_{n+1}^j$ for all $i$ and $j$. For fixed $n$, the map $\hat{\gamma} = \hat{\gamma}_n$ is given by a matrix $(\hat{\gamma}_{ji})$, where we can always assume that $\hat{\gamma}_{ji} \in \Zz_{>0}$ (considered as a map $\hat{H}_n^i = \Zz \to \Zz = \hat{H}_{n+1}^j$). Then $\gamma_n = \hat{\gamma} + \tau + t$ for homomorphisms $\tau: \: {\rm Tor}(H_n) \to {\rm Tor}(H_{n+1})$ and $t: \: \hat{H}_n \to {\rm Tor}(H_{n+1})$. Here we think of $\hat{H}_n$ as a subgroup (actually a direct summand) of $H_n$. As explained in \cite[\S~6]{GLN}, given a positive constant $\Gamma_n$ depending on $n$, we can always arrange that 
\begin{equation}
\label{e:gammaGamma}
  (\hat{\gamma}_n)_{ji} \geq \Gamma_n \ {\rm for} \ {\rm all} \ i \ {\rm and} \ j.
\end{equation}
Also, let $K_n$ be finitely generated abelian groups and $\chi_n: \: K_n \to K_{n+1}$ homomorphisms such that $K \cong \ilim \gekl{K_n; \chi_n}$.

Let $F_n = \bigoplus_i F_n^i$ be C*-algebras, where $F_n^{\bm{i}}$ is a homogeneous C*-algebra of the form $F_n^{\bm{i}} = P_n^{\bm{i}} M_{\infty}(C(Z_n^{\bm{i}})) P_n^{\bm{i}}$ for a connected compact space $Z_n^{\bm{i}}$ with base point $\theta_n^{\bm{i}}$ and a projection $P_n^{\bm{i}} \in M_{\infty}(C(Z_n^{\bm{i}}))$, while for all other indices $i \neq \bm{i}$, $F_n^i$ is a matrix algebra, $F_n^i = M_{[n,i]}$. We require that $(K_0(F_n^{\bm{i}}), K_0(F_n^{\bm{i}})^+,[1_{F_n^{\bm{i}}}]) \cong (H_n^{\bm{i}},(H_n^{\mathbf{i}})^+,([n,\bm{i}],\tau_n))$ and $K_1(F_n^{\bm{i}}) \cong K_n$, so that $(K_0(F_n), K_0(F_n)^+,[1_{F_n}], K_1(F_n)) \cong (H_n,H_n^+,u_n,K_n)$.

Let $\psi_n$ be a unital homomorphism $F_n \to F_{n+1}$ which induces $\gamma_n$ in $K_0$ and $\chi_n$ in $K_1$. We write $F_n = P_n C(Z_n) P_n$ where $Z_n = Z_n^{\bm{i}} \amalg \coprod_{i \neq \bm{i}} \{ \theta_n^i \}$, and $P_n = (P_n^{\bm{i}},(1_{[n,i]})_{i \neq \bm{i}}) \in M_{\infty}(C(Z_n^{\bm{i}})) \oplus \bigoplus_{i \neq \bm{i}} M_{[n,i]}(C(\{ \theta_n^i \}))$. Thus evaluation in $\theta_n^i$ induces a quotient map $\pi_n: \: F_n \to \hat{F}_n \defeq \bigoplus_i \hat{F}_n^i$, where $\hat{F}_n^i = M_{[n,i]}$. We require that $\psi_n$ induce homomorphisms $\hat{\psi}_n: \: \hat{F}_n \to \hat{F}_{n+1}$ so that we obtain a commutative diagram
\begin{equation*}
  \xymatrix{
  F_n \ar[d]_{\pi_n} \ar[r]_{\psi_n} & F_{n+1} \ar[d]^{\pi_{n+1}}
  \\
  \hat{F}_n \ar[r]^{\hat{\psi}_n} & \hat{F}_{n+1}
  }
\end{equation*}
which induces in $K_0$
\begin{equation*}
  \xymatrix{
  H_n \ar[d] \ar[r]_{\gamma_n} & H_{n+1} \ar[d]
  \\
  H_n / {\rm Tor}(H_n) \ar[r]^(.45){\hat{\gamma}_n} & H_{n+1} / {\rm Tor}(H_{n+1})
  }
\end{equation*}
where the vertical arrows are the canonical projections. As ${\rm Tor}(H_n) \subseteq G_n$, $H_n / G_n$ is torsion-free, and there is a canonical projection $H_n / {\rm Tor}(H_n) \to H_n / G_n$. Now let $E_n \defeq \bigoplus_p E_n^p$, $E_n^p = M_{\gekl{n,p}}$, so that $K_0(E_n) \cong H_n / G_n$, and for fixed $n$, let $\beta_0, \, \beta_1: \: \hat{F}_n \to E_n$ be unital homomorphisms which yield the commutative diagram
$$
  \xymatrix{
  K_0(\hat{F}_n) \ar[d]_{\cong} \ar[r]_{(\beta_0)_* - (\beta_1)_*} & K_0(E_n) \ar[d]^{\cong}
  \\
  H_n / {\rm Tor}(H_n) \ar[r]^(.55){\hat{\gamma}_n} & H_n / G_n.
  }
$$
We can assume $\beta_0 \oplus \beta_1: \: \hat{F}_n \to E_n \oplus E_n$ to be injective, because only the difference $(\beta_0)_* - (\beta_1)_*$ matters.

Define
\begin{eqnarray*}
  A_n &\defeq& \menge{(f,a) \in C([0,1],E_n) \oplus F_n}{f(t) = \beta_t(\pi_n(a)) \ {\rm for} \ t = 0,1},\\
  \hat{A}_n &\defeq& \menge{(f,\hat{a}) \in C([0,1],E_n) \oplus \hat{F}_n}{f(t) = \beta_t(\hat{a}) \ {\rm for} \ t = 0,1}.
\end{eqnarray*}
As $\beta_0 \oplus \beta_1$ is injective, we view $\hat{A}_n$ as a subalgebra of $C([0,1],E_n)$ via $(f,\hat{a}) \ma f$.

Choose for each $n$ a unital homomorphism $\hat{\varphi}_n: \: \hat{A}_n \to \hat{A}_{n+1}$ such that the composition with the map $C([0,1],E_{n+1}) \onto C([0,1],E_{n+1}^q)$ induced by the canonical projection $E_{n+1} \onto E_{n+1}^q$,
$$
  \hat{A}_n \overset{\hat{\varphi}_n}{\lori} \hat{A}_{n+1} \into C([0,1],E_{n+1}) \onto C([0,1],E_{n+1}^q),
$$
is of the form
\begin{equation}
\label{e:phihat}
  C([0,1],E_n) \supseteq \hat{A}_n \ni f \ma 
    u^*
    \begin{pmatrix}
    V(f) & \\
    & D(f)
    \end{pmatrix}
    u,
\end{equation}
where $u$ is a continuous path of unitaries $[0,1] \to U(E_{n+1}^q)$,
$$
  V(f) = 
    \begin{pmatrix}
    \pi_1(f) & & \\
    & \pi_2(f) & \\
    & & \ddots
    \end{pmatrix}
$$
for some $\pi_{\bullet}$ of the form $\pi_{\bullet}: \: \hat{A}_n \to \hat{F}_n \onto \hat{F}_n^i$, where the first map is given by $(f,\hat{a}) \ma \hat{a}$ and the second map is the canonical projection, and
$$
  D(f) = 
    \begin{pmatrix}
    f \circ \lambda_1 & & \\
    & f \circ \lambda_2 & \\
    & & \ddots
    \end{pmatrix}
$$
for some continuous maps $\lambda_{\bullet}: \: [0,1] \to [0,1]$ with $\lambda_{\bullet}^{-1}(\gekl{0,1}) \subseteq \gekl{0,1}$. We require that the diagram
$$
  \xymatrix{
  \hat{A}_n \ar[r]_{\hat{\varphi}_n} \ar[d] & \hat{A}_{n+1} \ar[d] \\
  \hat{F}_n \ar[r]^{\hat{\psi}_n} & \hat{F}_{n+1}
  }
$$
commute, where the vertical maps are given by $(f,\hat{a}) \ma \hat{a}$.

Then there exists a unique homomorphism $\varphi_n: \: A_n \to A_{n+1}$ which fits into the commutative diagram
$$
  \xymatrix{
  A_n \ar[rr]_{\varphi_n} \ar[dd] \ar[dr] & & A_{n+1} \ar[dr] \ar[dd] & \\
  & \hat{A}_n \ar[rr]^(.35){\hat{\varphi}_n} \ar[dd] & & \hat{A}_{n+1} \ar[dd] \\
  F_n \ar[rr]_(.6){\psi_n} \ar[dr] & & F_{n+1} \ar[dr] & \\
  & \hat{F}_n \ar[rr]^{\hat{\psi}_n} & & \hat{F}_{n+1}
  }
$$
where all the unlabelled arrows are given by the canonical maps.

By construction, $\ilim \gekl{A_n;\varphi_n}$ has the desired Elliott invariant (in particular, the canonical map $\ilim \gekl{A_n;\varphi_n} \to \hat{F} \defeq \ilim \gekl{\hat{F}_n; \hat{\psi}_n}$ induces $T(\ilim \gekl{A_n;\varphi_n}) \cong T(\hat{F})$). However, this is not a simple C*-algebra. Thus a further modification is needed to enforce simplicity. To this end, choose $\bm{I}_n \subseteq (0,1)$ and $\bm{Z}_n^{\bm{i}} \subseteq Z_n^{\bm{i}}$ $\frac{1}{n}$-dense and replace $\varphi_n: \: A_n \to A_{n+1}$ by the unital homomorphism $\xi_n: \: A_n \to A_{n+1}$ such that:
\setlength{\parindent}{0cm} \setlength{\parskip}{0cm}

\begin{itemize}
\item the compositions
$$
  A_n \overset{\xi_n}{\lori} A_{n+1} \to F_{n+1} \onto F_{n+1}^j\ \ \ {\rm and} \ \ \ A_n \overset{\varphi_n}{\lori} A_{n+1} \to F_{n+1} \onto F_{n+1}^j \ \text{coincide except for one index} \ j_{\xi} \neq \bm{j};
$$
\item the composition
$$
  A_n \overset{\xi_n}{\lori} A_{n+1} \to F_{n+1} \onto F_{n+1}^{j_{\xi}}  
$$
is of the form
$$
  A_n \ni (f,a) \ma 
    u^*
    \begin{pmatrix}
    I(f) & & \\
    & Z(a) & \\
    & & P(a)
    \end{pmatrix}
    u, 
$$
where $u$ is a permutation matrix in $M_{[n+1,j_{\xi}]}$,
$$
  I(f) = 
    \begin{pmatrix}
    f^{p_1}(t_1) & & \\
    & f^{p_2}(t_2) & \\
    & & \ddots
    \end{pmatrix}
$$
for indices $p_{\bullet}$ and $t_{\bullet} \in \bm{I}_n$ such that all possible pairs $p_{\bullet}, t_{\bullet}$ appear ($f^p$ is the component of $f$ in $C([0,1],E_n^p)$),
\begin{equation}
\label{e:Z(a)}
  Z(a) = 
    \begin{pmatrix}
    \tau_1(a(z_1)) & & \\
    & \tau_2(a(z_2)) & \\
    & & \ddots
    \end{pmatrix}
\end{equation}
for $z_{\bullet} \in \bm{Z}_n$ and isomorphisms $\tau_{\bullet}: \: P_n^{\bm{i}}(z_{\bullet})M_{\infty}P_n^{\bm{i}}(z_{\bullet}) \cong \hat{F}_n^{\bm{i}} = M_{[n,\bm{i}]}$, and
$$
  P(a) = 
    \begin{pmatrix}
    \pi_n^{i_1}(a) & & \\
    & \pi_n^{i_2}(a) & \\
    & & \ddots
    \end{pmatrix},
$$
where $\pi_n^i$ is the canonical projection $F_n \onto \hat{F}_n \onto \hat{F}_n^i$;
\item for every $q$, the composition
$$
  A_n \overset{\xi_n}{\lori} A_{n+1} \to C([0,1],E_{n+1}) \onto C([0,1],E_{n+1}^q)
$$
is of the form
$$
  A_n \ni (f,a) \ma 
    u^*
    \begin{pmatrix}
    \Phi(f) & \\
    & \Xi(a)
    \end{pmatrix}
    u,
$$
where $u$ is a continuous path of unitaries $[0,1] \to U(E_{n+1}^q)$, $\Phi(f)$ is of the same form
$$
  \begin{pmatrix}
  V(f) & \\
  & D(f)
  \end{pmatrix}
$$
as in \eqref{e:phihat}, 
$$
  \Xi(a)(t) = 
    \begin{pmatrix}
    \tau_1(t)(a(z_1(t))) & & \\
    & \tau_2(t)(a(z_2(t))) & \\
    & & \ddots
    \end{pmatrix}
$$
for continuous maps $z_{\bullet}: \: [0,1] \to Z_n^{\bm{i}}$, each of which is either a constant map with value in $\bm{Z}_n$ or connects $\theta_n^{\bm{i}}$ with $z_{\bullet} \in \bm{Z}_n$, and isomorphisms $\tau_{\bullet}(t): \: P_n^{\bm{i}}(z_{\bullet}(t)) M_{\infty}P_n^{\bm{i}}(z_{\bullet}(t)) \cong \hat{F}_n^{\bm{i}}$ depending continuously on $t \in [0,1]$ such that for $t \in \gekl{0,1}$, $\tau_{\bullet}(t) = \id$ if $z_{\bullet}(t) = \theta_n^{\bm{i}}$ and $\tau_{\bullet}(t) = \tau_{\bullet}$ if $z_{\bullet}(t) = z_{\bullet}$, where $\tau_{\bullet}$ is as in \eqref{e:Z(a)}.
\end{itemize}
\setlength{\parindent}{0cm} \setlength{\parskip}{0.5cm}

Then $\ilim \gekl{A_n; \xi_n}$ is a simple unital C*-algebra with prescribed Elliott invariant.

\subsection{The stably projectionless case}
\label{ss:StProjless}

We follow \cite{Ell} (see also \cite{EV}), with slight modifications as in the unital case. Let $(G_0,T,\rho,G_1)$ be as in Theorem~\ref{THM2}. 

Write $G = G_0$ and $K = G_1$. Choose a dense subgroup $G' \subseteq {\rm Aff}(T)$. Set $H \defeq G \oplus G'$, 
$$
  H^+ \defeq \gekl{(0,0)} \cup \menge{(g,f) \in G \oplus G'}{\rho(g)(\tau) + f(\tau) > 0 \ \forall \ \tau \in T}.
$$
Then $(H,H^+)$ becomes a simple ordered group, inducing the structure of a dimension group on $H / {\rm Tor}(H)$. Now construct a commutative diagram
$$
  \xymatrix{
  G_1 \ar[d] \ar[r] & \dotso \ar[r] & G_n \ar[d] \ar[r] & G_{n+1} \ar[d] \ar[r] & & \dotso & \ar[r] & G \ar[d]
  \\
  H_1 \ar[d] \ar[r]^{\gamma_1} & \dotso \ar[r]^{\gamma_{n-1}} & H_n \ar[d] \ar[r]^{\gamma_n} & H_{n+1} \ar[d] \ar[r]^{\gamma_{n+1}} & & \dotso & \ar[r] & H \ar[d]
  \\
  H_1 / G_1 \ar[r] & \dotso \ar[r] & H_n / G_n \ar[r] & H_{n+1} / G_{n+1} \ar[r] & & \dotso & \ar[r] & H / G
  }
$$
where
\setlength{\parindent}{0cm} \setlength{\parskip}{0cm}

\begin{itemize}
\item $H_n$ is a finitely generated abelian group with $H_n = \bigoplus_i H_n^i$, where for one distinguished index $\bm{i}$, $H_n^{\bm{i}} = \Zz \oplus {\rm Tor}(H_n)$, and for all other indices, $H_n^i = \Zz$;
\item with $(H_n^{\bm{i}})^+ \defeq \gekl{(0,0)} \cup (\Zz_{>0} \oplus {\rm Tor}(H_n))$, $(H_n^i)^+ \defeq \Zz_{\geq 0}$ for all $i \neq \bm{i}$, and $H_n^+ \defeq \bigoplus_i (H_n^i)^+ \subseteq H_n^{\bm{i}} \oplus \bigoplus_{i \neq \bm{i}} H_n^i = H_n$ we have 
\begin{equation}
\label{e:H=H_stprojless}
\ilim \gekl{(H_n,H_n^+); \gamma_n} \cong (H,H^+);
\end{equation}
\item with $G_n \defeq (\gamma_n^{\infty})^{-1}(G)$, where $\gamma_n^{\infty}: \: H_n \to H$ is the map provided by \eqref{e:H=H_stprojless}, we have $G_n \cap H_n^+ = \gekl{0}$, and \eqref{e:H=H} induces $\ilim \gekl{G_n; \gamma_n} \cong G$;
\item the vertical maps are the canonical ones.
\end{itemize}
\setlength{\parindent}{0cm} \setlength{\parskip}{0.5cm}

Let $\hat{\gamma}_n: \: H_n / {\rm Tor}(H_n) \eqdef \hat{H}_n = \bigoplus_i \hat{H}_n^i \to \bigoplus_j \hat{H}_{n+1}^j = \hat{H}_{n+1} \defeq H_{n+1} / {\rm Tor}(H_{n+1})$ be the homomorphism induced by $\gamma_n$, where $\hat{H}_n^i = \Zz = \hat{H}_{n+1}^j$ for all $i$ and $j$. For fixed $n$, the map $\hat{\gamma} = \hat{\gamma}_n$ is given by a matrix $(\hat{\gamma}_{ji})$, where we can always assume that $\hat{\gamma}_{ji} \in \Zz_{>0}$ (considered as a map $\hat{H}_n^i = \Zz \to \Zz = \hat{H}_{n+1}^j$). Then $\gamma_n = \hat{\gamma} + \tau + t$ for homomorphisms $\tau: \: {\rm Tor}(H_n) \to {\rm Tor}(H_{n+1})$ and $t: \: \hat{H}_n \to {\rm Tor}(H_{n+1})$. Here we think of $\hat{H}_n$ as a subgroup of $H_n$. As in the unital case (see \cite[\S~6]{GLN}), given a positive constant $\Gamma_n$ depending on $n$, we can always arrange that 
\begin{equation}
\label{e:gammaGamma_stprojless}
  (\hat{\gamma}_n)_{ji} \geq \Gamma_n \ {\rm for} \ {\rm all} \ i \ {\rm and} \ j.
\end{equation}
Also, let $K_n$ be finitely generated abelian groups and $\chi_n: \: K_n \to K_{n+1}$ homomorphisms such that $K \cong \ilim \gekl{K_n; \chi_n}$.

Let $F_n = \bigoplus_i F_n^i$ be C*-algebras, where $F_n^{\bm{i}}$ is a homogeneous C*-algebra of the form $F_n^{\bm{i}} = P_n^{\bm{i}} M_{\infty}(C(Z_n^{\bm{i}})) P_n^{\bm{i}}$ for a connected compact space $Z_n^{\bm{i}}$ with base point $\theta_n^{\bm{i}}$ and a projection $P_n^{\bm{i}} \in M_{\infty}(C(Z_n^{\bm{i}}))$, while for all other indices $i \neq \bm{i}$, $F_n^i$ is a matrix algebra, $F_n^i = M_{[n,i]}$. We require that $(K_0(F_n^{\bm{i}}), K_0(F_n^{\bm{i}})^+) \cong (H_n^{\bm{i}},(H_n^{\mathbf{i}})^+)$ and $K_1(F_n^{\bm{i}}) \cong K_n$, so that $(K_0(F_n), K_0(F_n)^+, K_1(F_n)) \cong (H_n,H_n^+,K_n)$.

Let $\psi_n$ be a unital homomorphism $F_n \to F_{n+1}$ which induces $\gamma_n$ in $K_0$ and $\chi_n$ in $K_1$. We write $F_n = P_n C(Z_n) P_n$ where $Z_n = Z_n^{\bm{i}} \amalg \coprod_{i \neq \bm{i}} \{ \theta_n^i \}$, and $P_n = (P_n^{\bm{i}},(1_{[n,i]})_{i \neq \bm{i}}) \in M_{\infty}(C(Z_n^{\bm{i}})) \oplus \bigoplus_{i \neq \bm{i}} M_{[n,i]}(C(\{ \theta_n^i \}))$. Thus, evaluation in $\theta_n^i$ induces a quotient map $\pi_n: \: F_n \to \hat{F}_n \defeq \bigoplus_i \hat{F}_n^i$, where $\hat{F}_n^i = M_{[n,i]}$. We require that $\psi_n$ induce homomorphisms $\hat{\psi}_n: \: \hat{F}_n \to \hat{F}_{n+1}$ so that we obtain a commutative diagram
\begin{equation*}
  \xymatrix{
  F_n \ar[d]_{\pi_n} \ar[r]_{\psi_n} & F_{n+1} \ar[d]^{\pi_{n+1}}
  \\
  \hat{F}_n \ar[r]^{\hat{\psi}_n} & \hat{F}_{n+1}
  }
\end{equation*}
which induces in $K_0$
\begin{equation*}
  \xymatrix{
  H_n \ar[d] \ar[r]_{\gamma_n} & H_{n+1} \ar[d]
  \\
  H_n / {\rm Tor}(H_n) \ar[r]^(.45){\hat{\gamma}_n} & H_{n+1} / {\rm Tor}(H_{n+1})
  }
\end{equation*}
where the vertical arrows are the canonical projections. As ${\rm Tor}(H_n) \subseteq G_n$, $H_n / G_n$ is torsion-free, and there is a canonical projection $H_n / {\rm Tor}(H_n) \to H_n / G_n$. Now let $E_n \defeq \bigoplus_p E_n^p$, $E_n^p = M_{\gekl{n,p}}$, such that $K_0(E_n) \cong H_n / G_n$, and for fixed $n$, let $\beta_0, \, \beta_1: \: \hat{F}_n \to E_n$ be (necessarily non-unital) homomorphisms which yield the commutative diagram
$$
  \xymatrix{
  K_0(\hat{F}_n) \ar[d]_{\cong} \ar[r]_{(\beta_0)_* - (\beta_1)_*} & K_0(E_n) \ar[d]^{\cong}
  \\
  H_n / {\rm Tor}(H_n) \ar[r]^(.55){\hat{\gamma}_n} & H_n / G_n.
  }
$$
As in the unital case, we may assume $\beta_0 \oplus \beta_1: \: \hat{F}_n \to E_n \oplus E_n$ to be injective.

Define
\begin{eqnarray*}
  A_n &\defeq& \menge{(f,a) \in C([0,1],E_n) \oplus F_n}{f(t) = \beta_t(\pi_n(a)) \ {\rm for} \ t = 0,1},\\
  \hat{A}_n &\defeq& \menge{(f,\hat{a}) \in C([0,1],E_n) \oplus \hat{F}_n}{f(t) = \beta_t(\hat{a}) \ {\rm for} \ t = 0,1}.
\end{eqnarray*}
As $\beta_0 \oplus \beta_1$ is injective, we view $\hat{A}_n$ as a subalgebra of $C([0,1],E_n)$ via $(f,\hat{a}) \ma f$.

Choose for each $n$ a homomorphism $\hat{\varphi}_n: \: \hat{A}_n \to \hat{A}_{n+1}$ such that the composition with the map $C([0,1],E_{n+1}) \onto C([0,1],E_{n+1}^q)$ induced by the canonical projection $E_{n+1} \onto E_{n+1}^q$,
$$
  \hat{A}_n \overset{\hat{\varphi}_n}{\lori} \hat{A}_{n+1} \into C([0,1],E_{n+1}) \onto C([0,1],E_{n+1}^q),
$$
is of the form
$$
  C([0,1],E_n) \supseteq \hat{A}_n \ni f \ma 
    u^*
    \begin{pmatrix}
    V(f) & \\
    & D(f)
    \end{pmatrix}
    u,
$$
where $u$ is a continuous path of unitaries $[0,1] \to U(E_{n+1}^q)$,
$$
  V(f) = 
    \begin{pmatrix}
    \pi_1(f) & & \\
    & \pi_2(f) & \\
    & & \ddots
    \end{pmatrix}
$$
for some $\pi_{\bullet}$ of the form $\pi_{\bullet}: \: \hat{A}_n \to \hat{F}_n \onto \hat{F}_n^i$, where the first map is given by $(f,\hat{a}) \ma \hat{a}$ and the second map is the canonical projection, and
$$
  D(f) = 
    \begin{pmatrix}
    f \circ \lambda_1 & & \\
    & f \circ \lambda_2 & \\
    & & \ddots
    \end{pmatrix}
$$
for some continuous maps $\lambda_{\bullet}: \: [0,1] \to [0,1]$ with $\lambda_{\bullet}^{-1}(\gekl{0,1}) \subseteq \gekl{0,1}$. We require that 
$$
  \xymatrix{
  \hat{A}_n \ar[r]_{\hat{\varphi}_n} \ar[d] & \hat{A}_{n+1} \ar[d] \\
  \hat{F}_n \ar[r]^{\hat{\psi}_n} & \hat{F}_{n+1}
  }
$$
commutes, where the vertical maps are given by $(f,\hat{a}) \ma \hat{a}$.

Then there exists a unique homomorphism $\varphi_n: \: A_n \to A_{n+1}$ which fits into the commutative diagram
$$
  \xymatrix{
  A_n \ar[rr]_{\varphi_n} \ar[dd] \ar[dr] & & A_{n+1} \ar[dr] \ar[dd] & \\
  & \hat{A}_n \ar[rr]^(.35){\hat{\varphi}_n} \ar[dd] & & \hat{A}_{n+1} \ar[dd] \\
  F_n \ar[rr]_(.6){\psi_n} \ar[dr] & & F_{n+1} \ar[dr] & \\
  & \hat{F}_n \ar[rr]^{\hat{\psi}_n} & & \hat{F}_{n+1}
  }
$$
where all the unlabelled arrows are given by the canonical maps.

By construction, $\ilim \gekl{A_n;\varphi_n}$ has the desired Elliott invariant (the details are as in the unital case, see \cite[\S~13]{GLN}). The same modification as in the unital case produces new connecting maps $\xi_n: \: A_n \to A_{n+1}$ such that $\ilim \gekl{A_n;\xi_n}$ is a simple (stably projectionless) C*-algebra with prescribed Elliott invariant. Moreover, choosing $\xi_n$ with the property that strictly positive elements are sent to strictly positive elements, $\ilim \gekl{A_n;\xi_n}$ will have continuous scale by \cite[Theorem~9.3]{GLI} (compare also \cite[\S~6]{GLII}). In addition, we choose $\xi_n$ such that full elements are sent to full elements.

\bremark
In an earlier version of this paper, we modified the construction in \cite[\S~6]{GLII} instead, which covers all Elliott invariants for stably projectionless C*-algebras with trivial pairing between K-theory and traces ($\rho = 0$). I would like to thank the referee for pointing out that \cite{Ell} (see also \cite{EV}) describes a general construction exhausting all possible Elliott invariants with weakly unperforated pairing between K-theory and traces (in the stably projectionless case, this is precisely the condition that $\rho$ is weakly unperforated as in Theorem~\ref{THM2}).
\eremark

\section{Concrete construction of AH-algebras}
\label{s:AH}

We start with the following standard fact.
\blemma
\label{L:XN}
Given an integer $N > 1$, let $\mu_N: \: S^1 \to S^1, \, z \ma z^N$, and set $X_N \defeq D^2 \cup_{\mu_N} S^1$, where we identify $z \in S^1 = \partial D^2$ with $\mu_N(z) \in S^1$. Then
$$
  H^{\bullet}(X_N) \cong 
  \bfa
  \Zz & {\rm if} \ \bullet = 0;\\
  \Zz / N & {\rm if} \ \bullet = 2;\\
  \gekl{0} & {\rm else}.
  \efa
$$
Moreover, $(K_0(C(X_N)),K_0(C(X_N))^+,[1_{C(X_N)}],K_1(C(X_N))) \cong (\Zz \oplus \Zz / N, \gekl{(0,0)} \cup (\Zz_{>0} \oplus \Zz / N), (1,0), \gekl{0})$.
\elemma

In the following, we view $S^2$ as the one point compactification of $\mathring{D}^2$, $S^2 = \mathring{D}^2 \cup \gekl{\infty}$.
\blemma
\label{L:pXN}
Let $X_N \onto S^2$ be the continuous map sending $\mathring{D}^2 \subseteq D^2$ identically to $\mathring{D}^2 \subseteq S^2$, $\partial D^2$ to $\infty$ and $S^1$ to $\infty$. Let $p_{X_N}$ be the pullback of the Bott line bundle on $S^2$ (see for instance \cite[\S~6.2]{Ros}) to $X_N$ via this map. We view $p_{X_N}$ as a projection in $M_2(C(X_N))$. Then there is an isomorphism $K_0(C(X_N)) \cong \Zz \oplus \Zz / N$ identifying the class of $1_{C(X_N)}$ with the generator of $\Zz$ and the class of $p_{X_N}$ with $(1,1)$.
\elemma
\setlength{\parindent}{0cm} \setlength{\parskip}{0cm}

\bproof
Just analyse the K-theory exact sequence attached to $0 \to C_0(\mathring{D}^2) \to C(X_N) \to C(S^1) \to 0$.
\eproof
\setlength{\parindent}{0cm} \setlength{\parskip}{0.5cm}

We recall another standard fact.
\blemma
\label{L:YN}
Given an integer $N > 1$, let $Y_N \defeq \Sigma X_N \cong D^3 \cup_{\Sigma \mu_N} S^2$, where we identify $z \in S^2 = \partial D^3 \cong \Sigma S^1$ with $(\Sigma \mu_N)(z) \in \Sigma S^1 \cong S^2$. (Here $\Sigma$ stands for suspension.) Then
$$
  H^{\bullet}(Y_N) \cong 
  \bfa
  \Zz & {\rm if} \ \bullet = 0;\\
  \Zz / N & {\rm if} \ \bullet = 3;\\
  \gekl{0} & {\rm else}.
  \efa
$$
Moreover, $K_0(C(Y_N)) = \Zz[1_{C(Y_N)}]$ and $K_1(C(Y_N)) \cong \Zz / N$.
\elemma

In the following, we view $S^3$ as the one point compactification of $\mathring{D}^3$, $S^3 = \mathring{D}^3 \cup \gekl{\infty}$.
\blemma
\label{L:Y-->S}
Let $Y_N \onto S^3$ be the continuous map sending $\mathring{D}^3 \subseteq D^3$ identically to $\mathring{D}^3 \subseteq S^3$, $\partial D^3$ to $\infty$ and $S^2$ to $\infty$. Then the dual map $C(S^3) \to C(Y_N)$ induces in $K_1$ a surjection $K_1(C(S^3)) \cong \Zz \onto \Zz / N \cong K_1(C(Y_N))$.
\elemma
\setlength{\parindent}{0cm} \setlength{\parskip}{0cm}

\bproof
Just analyse the K-theory exact sequence attached to $0 \to C_0(\mathring{D}^3) \to C(Y_N) \to C(S^2) \to 0$.
\eproof
\setlength{\parindent}{0cm} \setlength{\parskip}{0.5cm}

Analysing K-theory exact sequences, the following is a straightforward observation.
\blemma
\label{L:Psi_m}
Let $N, N' \in \Zz_{>1}$ and $m \in \Zz_{>0}$ with $N' \mid m \cdot N$, say $m \cdot N = m' \cdot N'$. Define a continuous map 
$$
  \Psi_m^*: \: X_{N'} = D^2 \cup_{\mu_{N'}} S^1 \to D^2 \cup_{\mu_N} S^1 = X_N
$$
by sending $x \in D^2$ to $x^m \in D^2$ and $z \in S^1$ to $z^{m'} \in S^1$. Then the dual map $\Psi_m: \: C(X_N) \to C(X_{N'})$ induces in $K_0$ the homomorphism
$$
  K_0(C(X_N)) \cong \Zz \oplus \Zz / N 
  \overset{
  \rukl{
  \begin{smallmatrix}
  1 & 0 \\
  0 & m
  \end{smallmatrix}
  }
  }{\lori} \Zz \oplus \Zz / N' \cong K_0(C(X_{N'})).
$$
\elemma

Naturality of suspension yields
\blemma
\label{L:SigmaPsi_m}
Let $N, N' \in \Zz_{>1}$ and $m \in \Zz_{>0}$ with $N' \mid m \cdot N$, say $m \cdot N = m' \cdot N'$. Let $\Sigma \Psi_m: \: C(Y_N) \to C(Y_{N'})$ be the map dual to $\Sigma \Psi_m^*: \: Y_{N'} \cong \Sigma X_{N'} \to \Sigma X_N \cong Y_N$. Then $\Sigma \Psi_m$ induces in $K_1$ the homomorphism
$$
  K_1(C(Y_N)) \cong \Zz / N \overset{m}{\lori} \Zz / N' \cong K_1(C(Y_{N'})).
$$
\elemma

In the following, we view $X_N$ and $Y_N$ as pointed spaces, with base point $1 = (1,0) \in S^1 = \partial D^2 \subseteq D^2$ in $X_N$ and base point $(1,0,0) \in S^2 = \partial D^3 \subseteq D^3$ in $Y_N$. Note that $\Psi_m$ and $\Sigma \Psi_m$ preserve base points. Moreover, if $\theta$ denotes the base point of $X_N$, then the projection $p_{X_N}$ in Lemma~\ref{L:pXN} satisfies
\begin{equation}
\label{e:pXN=_}
  p_{X_N}(\theta) = 
  \rukl{
  \begin{smallmatrix}
  1 & 0 \\
  0 & 0
  \end{smallmatrix}
  }.
\end{equation}

Now let $H_n = H_n^{\bm{i}} \oplus \bigoplus_{i \neq \bm{i}} H_n^i$, $H_{n+1} = H_{n+1}^{\bm{j}} \oplus \bigoplus_{j \neq \bm{j}} H_{n+1}^j$ be abelian groups with $H_n^{\bm{i}} = \Zz \oplus T_n$, $H_{n+1}^{\bm{j}} = \Zz \oplus T_{n+1}$ for finitely generated torsion groups $T_n$, $T_{n+1}$, and $H_n^i = \Zz$, $H_{n+1}^j = \Zz$ for all $i \neq \bm{i}$, $j \neq \bm{j}$. Let $(H_n^{\bm{i}})^+ \defeq \gekl{(0,0)} \cup (\Zz_{>0} \oplus T_n)$, $(H_n^i)^+ \defeq \Zz_{\geq 0}$ for all $i \neq \bm{i}$, $H_n^+ \defeq \bigoplus_i (H_n^i)^+ \subseteq H_n^{\bm{i}} \oplus \bigoplus_{i \neq \bm{i}} H_n^i = H_n$ and $u_n = (([n,\bm{i}],\tau_n),([n,i])_{i \neq \bm{i}}) \in H_n^+$. Similarly, define $(H_{n+1}^j)^+$, $H_{n+1}^+ \defeq \bigoplus_i (H_{n+1}^j)^+$, and let $u_{n+1} = (([n+1,\bm{j}],\tau_{n+1}),([n+1,j])_{j \neq \bm{j}}) \in H_{n+1}^+$. Let $T_n = \bigoplus_k T_n^k$, where $T_n^k = \Zz / N_n^k$, and $T_{n+1} = \bigoplus_k T_{n+1}^l$, where $T_{n+1}^l = \Zz / N_{n+1}^l$. Let $\gamma_n: \: H_n \to H_{n+1}$ be a homomorphism with $\gamma_n(u_n) = u_{n+1}$. (In the stably projectionless case, these order units are not part of the given data, but we can always choose such order units.) Let us fix $n$, and suppose that $\gamma = \gamma_n$ induces a homomorphism $\hat{\gamma}: \: H_n / {\rm Tor}(H_n) = \hat{H}_n = \bigoplus_i \hat{H}_n^i \to \bigoplus_j \hat{H}_{n+1}^j = \hat{H}_{n+1} = H_{n+1} / {\rm Tor}(H_{n+1})$, where $\hat{H}_n^i = \Zz = \hat{H}_{n+1}^j$ for all $i$ and $j$. Viewing $\hat{H}_n$ as a subgroup (actually a direct summand) of $H_n$, we obtain that $\gamma_n = \hat{\gamma} + \tau + t$ for homomorphisms $\tau: \: {\rm Tor}(H_n) \to {\rm Tor}(H_{n+1})$ and $t: \: \hat{H}_n \to {\rm Tor}(H_{n+1})$. $\hat{\gamma}$ is given by an integer matrix $(\hat{\gamma}_{ji})$. Similarly, $\tau$ is given by an integer matrix $(\tau_{lk})$, where we view $\tau_{lk}$ as a homomorphism $T_n^k \to T_{n+1}^l$. Also, $t$ is given by an integer matrix $(t_{li})$, where we view $t_{li}$ as a homomorphism $H_n^i \to T_{n+1}^l$. Clearly, we can always arrange $\tau_{lk}, t_{li} > 0$ for all $l, k, i$, and because of \eqref{e:gammaGamma} and \eqref{e:gammaGamma_stprojless}, we can also arrange
\begin{equation}
\label{e:gammageq_}
  \hat{\gamma}_{ji} > 0 \ {\rm and} \ \hat{\gamma}_{\bm{j} \bm{i}} \geq \#_0(k) + 1.
\end{equation}
Here $\#_0(k)$ is the number of direct summands in $T_n$ (i.e., the number of indices $k$).

We have the following direct consequence of Lemma~\ref{L:XN}.
\blemma
\label{L:KXn}
Let $X_n^{\bm{i}} \defeq \bigvee_k X_{N_n^k}$, where we take the wedge sum with respect to the base points of the individual $X_{N_n^k}$. Denote the base point of $X_n^{\bm{i}}$ by $\theta_n^{\bm{i}}$. Set $X_n \defeq X_n^{\bm{i}} \amalg \coprod_{i \neq \bm{i}} \{ \theta_n^i \}$. Then
$$
  (K_0(C(X_n)),K_0(C(X_n))^+,K_1(C(X_n))) \cong (H_n,H_n^+,\gekl{0}).
$$
\elemma
Define $X_{n+1}$ in an analogous way, i.e., $X_{n+1}^{\bm{j}} \defeq \bigvee_l X_{N_{n+1}^l}$, and $X_{n+1} \defeq X_{n+1}^{\bm{j}} \amalg \coprod_{j \neq \bm{j}} \{ \theta_{n+1}^j \}$. Now, for fixed $n$, our goal is to construct a homomorphism $\psi$ realizing the homomorphism $\gamma$ in $K_0$.

The map $\bigvee_l \Psi_{\tau_{lk}}^*: \: \bigvee_l X_{N_{n+1}^l} \to X_{N_n^k}$ induces the dual homomorphism $\psi_{\tau}^k: \: C(X_{N_n^k}) \to C(X_{n+1}^{\bm{j}})$. Here $\Psi_{\tau_{lk}}$ are the maps from Lemma~\ref{L:Psi_m}. The direct sum $\bigoplus_k \psi_{\tau}^k: \: \bigoplus_k C(X_{N_n^k}) \to M_{\#_0(k)}(C(X_{n+1}^{\bm{j}}))$ restricts to a homomorphism $\psi_{\tau}: \: C(X_n^{\bm{i}}) = C(\bigvee_k X_{N_n^k}) \to M_{\#_0(k)}(C(X_{n+1}^{\bm{j}}))$.

Let $p^{(\bm{i})} \in M_2(C(X_{n+1}^{\bm{j}})) = M_2(C(\bigvee_l X_{N_{n+1}^l}))$ be given by $p^{(\bm{i})} \vert_{C(X_{N_{n+1}^l})} = M_2(\Psi_{t_{l \bm{i}}})(p_{X_{N_{n+1}^l}})$. Define $\psi_t$ as the composite
$$
  C(X_n^{\bm{i}}) \overset{\ev_{\theta_n^{\bm{i}}}}{\lori} \Cz \to M_2(C(X_{n+1}^{\bm{j}})), \ \text{where the second map is given by} \ 1 \ma p^{(\bm{i})}.
$$

Moreover, define $\psi_{\bm{j} \bm{i}}: \: C(X_n^{\bm{i}}) \to M_{\hat{\gamma}_{\bm{j} \bm{i}} + 1}(C(X_{n+1}^{\bm{j}}))$ by setting
$$
  \psi_{\bm{j} \bm{i}}(f) = 
    \begin{pmatrix}
    f(\theta_n^{\bm{i}}) & & & & \\
    & \ddots & & & \\
    & & f(\theta_n^{\bm{i}}) & & \\
    & & & \psi_{\tau}(f) & \\
    & & & & \psi_t(f) 
    \end{pmatrix}
$$
where we put $\hat{\gamma}_{\bm{j} \bm{i}} - \#_0(k) - 1$ copies of $f(\theta_n^{\bm{i}})$ on the diagonal.

For $i \neq \bm{i}$, let $p^{(i)} \in M_2(C(X_{n+1}^{\bm{j}})) = M_2(C(\bigvee_l X_{N_{n+1}^l}))$ be given by $p^{(i)} \vert_{C(X_{N_{n+1}^l})} = M_2(\Psi_{t_{l i}})(p_{X_{N_{n+1}^l}})$. Define
$$
  \psi_{\bm{j} i}: \: C(\{ \theta_n^i \}) = \Cz \to M_{\hat{\gamma}_{\bm{j} i} + 1}(C(X_{n+1}^{\bm{j}}))
\ \text{by sending} \ 1 \in \Cz \ {\rm to} \
    \begin{pmatrix}
    1 & & & \\
    & \ddots & & \\
    & & 1 & \\
    & & & p^{(i)} 
    \end{pmatrix},
$$
where we put $\hat{\gamma}_{\bm{j} i} - 1$ copies of $1$ on the diagonal.

For $j \neq \bm{j}$, define 
$$
  \psi_{j \bm{i}}: \: C(X_n^{\bm{i}}) \to M_{\hat{\gamma}_{j \bm{i}}}(C(\{ \theta_{n+1}^{j} \})), \, f \ma 
    \begin{pmatrix}
    f(\theta_n^{\bm{i}}) & & \\
    & \ddots & \\
    & & f(\theta_n^{\bm{i}})
    \end{pmatrix},  
$$
where we put $\hat{\gamma}_{j \bm{i}}$ copies of $f(\theta_n^{\bm{i}})$ on the diagonal.

For $i \neq \bm{i}$ and $j \neq \bm{j}$, define
$$
  \psi_{ji}: \: C(\{ \theta_n^{i} \}) \to M_{\hat{\gamma}_{ji}}(C(\{ \theta_{n+1}^{j} \})), \, 1 \ma 
    \begin{pmatrix}
    1 & & \\
    & \ddots & \\
    & & 1
    \end{pmatrix},  
$$
where we put $\hat{\gamma}_{ji}$ copies of $1$ on the diagonal.

To unify notation, let us set $X_n^i = \{ \theta_n^i \}$, $X_{n+1}^j = \{ \theta_{n+1}^j \}$.

For $u_n = (([n,\bm{i}],\tau_n),([n,i])_{i \neq \bm{i}}) \in H_n^+$, let $s(n,\bm{i})$ be a positive integer and $P_n^{\bm{i}} \in M_{s(n,\bm{i})}(C(X_n^{\bm{i}}))$ a projection such that:
\setlength{\parindent}{0cm} \setlength{\parskip}{0cm}

\begin{itemize}
\item $P_n^{\bm{i}}$ is a sum of line bundles;
\item $[P_n^{\bm{i}}]$ corresponds to $([n,\bm{i}],\tau_n)$ under the identification in Lemma~\ref{L:KXn};
\item $P_n^{\bm{i}}(\theta_n^{\bm{i}}) = 1_{[n,\bm{i}]}$ is of the form
$$
  u^*
    \begin{pmatrix}
    1 & & & & & \\
    & \ddots & & & & \\
    & & 1 & & & \\
    & & & 0 & & \\
    & & & & \ddots & \\
    & & & & & 0
    \end{pmatrix}
  u, 
  \ \ \ {\rm where} \ u \ \text{is a permutation matrix}.
$$
\end{itemize}
\setlength{\parindent}{0cm} \setlength{\parskip}{0.5cm}

$P_n^{\bm{i}}$ exists because of Lemma~\ref{L:pXN}. Moreover, we can extend $P_n^{\bm{i}}$ by $1_{[n,i]}$ to a projection in $\bigoplus_i M_{s(n,i)}(C(X_n^i))$ such that $[P_n]$ corresponds to $u_n$ under the isomorphism in Lemma~\ref{L:KXn}. Here $s(n,i) = [n,i]$ whenever $i \neq \bm{i}$. Then
\begin{equation}
\label{e:M(psi)}
  \rukl{M_{s(n,i)}(\psi_{ji})}_{ji}: \: \bigoplus_i M_{s(n,i)}(C(X_n^i)) \to \bigoplus_j M_{s(n+1,j)}(C(X_{n+1}^j))
\end{equation}
sends $P_n$ to $P_{n+1}$, where $P_{n+1}$ is of the same form as $P_n$, with $[P_{n+1}]$ corresponding to $u_{n+1}$ under the isomorphism from Lemma~\ref{L:KXn}. Hence the map in \eqref{e:M(psi)} restricts to a unital homomorphism
\begin{equation}
\label{e:M(psi)-res}
  P_n M_{\infty}(C(X_n)) P_n \to P_{n+1} M_{\infty}(C(X_{n+1})) P_{n+1}
\end{equation}
which in $K_0$ induces $\gamma$ by Lemma~\ref{L:Psi_m}. 

Now we turn to $K_1$. Assume $K_n = \bigoplus_i K_n^i$ is an abelian group, where for a distinguished index $\bm{i}$, $K_n^{\bm{i}} = T_n$ is a finitely generated torsion group $T_n = \bigoplus_k T_n^k$, $T_n^k = \Zz / N_n^k$, and $K_n^i = \Zz$ for all $i \neq \bm{i}$. Similarly, let $K_{n+1} = \bigoplus_j K_{n+1}^j$ be an abelian group, where for a distinguished index $\bm{j}$, $K_{n+1}^{\bm{j}} = T_{n+1}$ is a finitely generated torsion group $T_{n+1} = \bigoplus_l T_{n+1}^l$, $T_{n+1}^l = \Zz / N_{n+1}^l$, and $K_{n+1}^j = \Zz$ for all $j \neq \bm{j}$. For fixed $n$, let $\chi: \: K_n \to K_{n+1}$ be a homomorphism which is a sum $\chi = \hat{\chi} + \tau + t$, where $\hat{\chi}: \: \bigoplus_{i \neq \bm{i}} K_n^i \to \bigoplus_{j \neq \bm{j}} K_{n+1}^j$ is given by an integer matrix $(\hat{\chi}_{ji})$ (viewing $\hat{\chi}_{ji}$ as a homomorphism $K_n^i \to K_{n+1}^j$), $\tau: \: T_n \to T_{n+1}$ is given by an integer matrix $(\tau_{lk})$ (viewing $\tau_{lk}$ as a homomorphism $T_n^k \to T_{n+1}^l$), and $t: \: \bigoplus_{i \neq \bm{i}} K_n^i \to T_{n+1}$ is given by an integer matrix $(t_{li})$ (viewing $t_{li}$ as a homomorphism $K_n^i \to T_{n+1}^l$). We can always arrange that all the entries of these matrices are positive integers.

The following is a direct consequence of Lemma~\ref{L:YN}.
\blemma
Let $Y_n^{\bm{i}} = \bigvee_k Y_{N_n^k}$ and $Y_n = Y_n^{\bm{i}} \vee \bigvee_{i \neq \bm{i}} S^3$. Then $K_0(C(Y_n)) \cong \Zz$ and $K_1(C(Y_n)) \cong K_n$.
\elemma
We view $Y_n$ as a pointed space, and let $\theta_n$ be the base point of $Y_n$. Now let $\psi_{\tau}^k: \: C(Y_{N_n^k}) \to C(\bigvee_l Y_{N_{n+1}^l}) = C(Y_{n+1}^{\bm{j}})$ be the dual homomorphism of the map $\bigvee_l \Sigma (\Psi_{\tau_{lk}}^*): \: Y_{n+1}^{\bm{j}} = \bigvee_l Y_{N_{n+1}^l} \to Y_{N_n^k}$. Here $\Sigma (\Psi_{\tau_{lk}}^*)$ are the maps from Lemma~\ref{L:SigmaPsi_m}. The direct sum $\bigoplus_k \psi_{\tau}^k: \: \bigoplus_k C(Y_{N_n^k}) \to M_{\#_1(k)}(C(Y_{n+1}^{\bm{j}}))$ restricts to a homomorphism $\psi_{\bm{i}}: \: C(Y_n^{\bm{i}}) = C(\bigvee_k Y_{N_n^k}) \to M_{\#_1(k)}(C(Y_{n+1}^{\bm{j}})) \into M_{\#_1(k)}(C(Y_{n+1}))$.

For $i \neq \bm{i}$, define $\psi_{\bm{j}i}: \: C(Y_n^i) = C(S^3) \to C(Y_{n+1}^{\bm{j}})$ as the dual map of the composite
$$
  Y_{n+1}^{\bm{j}} = \bigvee_l Y_{N_{n+1}^l} \overset{\bigvee_l \Sigma(\psi_{t_{li}}^*)}{\lori} \bigvee_l Y_{N_{n+1}^l} \overset{\bigvee_l \Omega_l^*}{\lori} S^3,
$$
where $\Omega^*_l$ is the map $Y_{N_{n+1}^l} \to S^3$ constructed in Lemma~\ref{L:Y-->S}.

For $i \neq \bm{i}$ and $j \neq \bm{j}$, define $\psi_{ji}: \: C(Y_n^i) = C(S^3) \to C(S^3) = C(Y_{n+1}^j)$ as the dual map of $\Sigma \Sigma \mu_{\hat{\chi}_{ji}}: \: S^3 \cong \Sigma \Sigma S^1 \to \Sigma \Sigma S^1 \cong S^3$, where $\mu_{\hat{\chi}_{ji}}$ is the map from Lemma~\ref{L:XN}.

For every $i \neq \bm{i}$, we thus obtain the direct sum $\bigoplus_j \psi_{ji}: \: C(Y_n^i) \to \bigoplus_j C(Y_{n+1}^j)$ with image in $C(Y_{n+1}) = C(\bigvee_j Y_{n+1}^j) \subseteq \bigoplus_j C(Y_{n+1}^j)$. Hence we obtain a homomorphism $\psi_i: \: C(Y_n^i) \to C(Y_{n+1})$.

Now let $\#_1(i)$ be the number of summands of $K_n$. Then let $\psi: \: C(Y_n) \to M_{\#_1(k) + \#_1(i) - 1}(C(Y_{n+1}))$ be the restriction of $\bigoplus_i \psi_i$ to $C(Y_n) = C(\bigvee_i Y_n^i) \subseteq \bigoplus_i C(Y_n^i)$. By construction, and using Lemmas~\ref{L:Y-->S} and \ref{L:SigmaPsi_m}, $\psi$ induces $\chi$ in $K_1$.

We now combine our two constructions. Define $Z_n = X_n \vee Y_n$, where we identify the base point $\theta_n^{\bm{i}} \in X_n^{\bm{i}} \subseteq X_n$ with $\theta_n \in Y_n$. We extend $P_n$ from $X_n$ constantly to $Y_n$ (with value $P_n(\theta_n^{\bm{i}})$). Note that $\rk(P_{n+1}(\theta_{n+1}^{\bm{j}})) = \hat{\gamma}_{\bm{j} \bm{i}} \cdot \rk(P_n (\theta_n^{\bm{i}}))$. Because of \eqref{e:gammaGamma} and \eqref{e:gammaGamma_stprojless}, we can arrange $\hat{\gamma}_{\bm{j} \bm{i}} \geq \#_1(k) + \#_1(i) - 1$. By adding $\ev_{\theta_n}$ on the diagonal if necessary, we can modify $\psi$ to a homomorphism $\psi: \: C(Y_n) \to M_{\hat{\gamma}_{\bm{j} \bm{i}}}(C(Y_{n+1}))$ which induces $\gamma$ in $K_1$. We can thus think of $M_{[n,\bm{i}]}(\psi)$ as a unital homomorphism $P_n(\theta_n^{\bm{i}}) M_{s(n,\bm{i})}(C(Y_n)) P_n(\theta_n^{\bm{i}}) \to P_{n+1}(\theta_{n+1}^{\bm{j}}) M_{s(n+1,\bm{j})}(C(Y_{n+1})) P_{n+1}(\theta_{n+1}^{\bm{j}})$, i.e., as a unital homomorphism $P_n M_{s(n,\bm{i})}(C(Y_n)) P_n \to P_{n+1} M_{s(n+1,\bm{j})}(C(Y_{n+1})) P_{n+1}$. In combination with the homomorphism \eqref{e:M(psi)-res}, we obtain a unital homomorphism   $$
  P_n M_{\infty}(C(Z_n)) P_n \to P_{n+1} M_{\infty}(C(Z_{n+1})) P_{n+1}
$$
which induces $\gamma$ in $K_0$, sending $u_n$ to $u_{n+1}$, and $\chi$ in $K_1$.

Evaluation at $\theta_n^{\bm{i}} = \theta_n$ and $\theta_n^i$ (for $i \neq \bm{i}$) induces a quotient homomorphism which fits into a commutative diagram
\begin{equation}
\label{e:psi-ev}
  \xymatrix{
  P_n M_{\infty}(C(Z_n)) P_n \ar[d] \ar[r] & P_{n+1} M_{\infty}(C(Z_{n+1})) P_{n+1} \ar[d]
  \\
  \bigoplus_i M_{[n,i]} \ar[r] & \bigoplus_j M_{[n+1,j]}
  }
\end{equation}
which induces in $K_0$
$$
  \xymatrix{
  H_n \ar[d] \ar[r]_{\gamma} & H_{n+1} \ar[d]
  \\
  \hat{H}_n \ar[r]^{\hat{\gamma}} & \hat{H}_{n+1}.
  }
$$

\bremark
\label{R:S3S1}
If all $K_n$ are torsion-free, then we can replace $S^3$ by $S^1$ in our construction of $Y_n$.
\eremark

\section{The complete construction}

\subsection{The general construction with concrete models}

Applying our construction in \S~\ref{s:AH}, we obtain concrete models for $F_n$, $\hat{F}_n$, $\gamma_n$ and $\hat{\gamma}_n$ which we now plug into the general construction in \S~\ref{ss:UnitalCase} and \S~\ref{ss:StProjless}. Note that it is crucial that we work with these concrete models from \S~\ref{s:AH}. The reason is that only for these models can we provide groupoid descriptions of the C*-algebras and their homomorphisms which arise in the general construction (see \S~\ref{s:GPDModels}).

Note that with these concrete models, the composition
$$
  M_{[n,i]} \into \hat{F}_n \overset{\beta_{\bullet}}{\lori} E_n \onto M_{\gekl{n,p}},
$$
where the first and third maps are the canonical ones, is of the form
$$
  x \ma 
  u^*
    \begin{pmatrix}
    x & & & & & \\
    & \ddots & & & & \\
    & & x & & & \\
    & & & 0 & & \\
    & & & & \ddots & \\
    & & & & & 0
    \end{pmatrix}
  u
$$
for a permutation matrix $u$.

Apart from inserting these concrete models, we keep the same construction as in \S~\ref{ss:UnitalCase} and \S~\ref{ss:StProjless}.

\subsection{Summary of the construction}
\label{ss:SummConst}

In both the unital and stably projectionless cases, the C*-algebra with the prescribed Elliott invariant which we constructed is an inductive limit whose building blocks are of the form
\begin{equation}
\label{e:A}
  A_n = \menge{(f,a) \in C([0,1],E_n) \oplus F_n}{f(t) = \beta_t(a) \ {\rm for} \ t = 0,1},
\end{equation}
where:
\setlength{\parindent}{0cm} \setlength{\parskip}{0cm}

\begin{itemize}
\item $E_n$ is finite dimensional; 
\item $F_n$ is homogeneous of the form $P_n M_{\infty}(Z_n) P_n$, where $P_n$ is a sum of line bundles, and there are points $\theta_n^i \in Z_n$, one for each connected component, and all connected components just consist of $\theta_n^i$ with the only possible exception being the component of a distinguished point $\theta_n^{\bm{i}}$;
\item both $\beta_0$ and $\beta_1$ are compositions of the form $F_n \to \bigoplus_i M_{[n,i]} \to E_n$, where the first homomorphism is given by evaluation in $\theta_n^i \in Z_n$ and the second homomorphism is determined by the composites $M_{[n,i]} \into \bigoplus_i M_{[n,i]} \to E_n \onto E_n^p$ (where $E_n^p$ is a matrix block of $E_n$), which are of the form
$$
  x \ma 
  v^*
    \begin{pmatrix}
    x & & & & & \\
    & \ddots & & & & \\
    & & x & & & \\
    & & & 0 & & \\
    & & & & \ddots & \\
    & & & & & 0
    \end{pmatrix}
  v
$$
for a permutation matrix $v$.
\end{itemize}
\setlength{\parindent}{0cm} \setlength{\parskip}{0.5cm}

The connecting maps $\varphi_n$ of our inductive limit can be described as two parts:
\begin{align}
\label{e:ConnMapPartOne}
  &A_n \to A_{n+1} \onto F_{n+1};\\
\label{e:ConnMapPartTwo}
  &A_n \to A_{n+1} \onto C([0,1],E_{n+1}).
\end{align}
Both parts send $(f,a) \in A_n$ to an element which is in diagonal form up to permutation, i.e.,
\begin{equation}
\label{e:FormOfConnMaps}
    u^*
    \begin{pmatrix}
    * & & \\
    & * & \\
    & & \ddots
    \end{pmatrix}
    u,
\end{equation}
where for the entries on the diagonal, there are the following possibilities:
\setlength{\parindent}{0cm} \setlength{\parskip}{0cm}

\begin{itemize}
\item a map of the form
\begin{equation}
\label{e:phi1}
  [0,1] \ni t \ma f^p(\lambda(t)), \ {\rm for} \ {\rm a} \ {\rm continuous} \ {\rm map} \ \lambda: \: [0,1] \to [0,1] \ {\rm with} \ \lambda^{-1}(\gekl{0,1}) \subseteq \gekl{0,1},
\end{equation}
where $f^p$ is the image of $f$ under the canonical projection $C([0,1],E_n) \onto C([0,1],E_n^p)$;
\item a map of the form
\begin{equation}
\label{e:phi2} 
  [0,1] \ni t \ma \tau(t) a(x(t)),
\end{equation}
where $x: \: [0,1] \to Z_n$ is continuous and $\tau(t): \: P_n(x(t)) M_{\infty} P_n(x(t)) \cong P_n(\theta_n^i) M_{\infty}P_n(\theta_n^i)$ is an isomorphism depending continuously on $t$, with $\theta_n^i$ in the same connected component as $x(t)$, and $\tau(t) = \id$ if $x(t) = \theta_n^i$;
\item an element of $P_{n+1} M_{\infty}(C(Z_{n+1})) P_{n+1}$ with support in an isolated point $\theta_{n+1}^j$, which is of the form
\begin{equation}
\label{e:phi3}
  f^p(\bm{t}), \ {\rm for} \ {\rm some} \ \bm{t} \in (0,1),
\end{equation}
where $f^p$ is the image of $f$ under the canonical projection $C([0,1],E_n) \onto C([0,1],E_n^p)$;
\item an element of $P_{n+1} M_{\infty}(C(Z_{n+1})) P_{n+1}$ with support in an isolated point $\theta_{n+1}^j$, which is of the form
\begin{equation}
\label{e:phi4}
  \tau(a(x)) \ \ \ \text{for some} \ x \in Z_n \ {\rm with} \ x \notin \{ \theta_n^i \}_i \ \text{and an isomorphism} \ \tau: \: P_n(x) M_{\infty} P_n(x) \cong P_n(\theta_n^i) M_{\infty} P_n(\theta_n^i),
\end{equation}
where $\theta_n^i$ is in the same connected component as $x$;
\item an element of $P_{n+1} M_{\infty}(C(Z_{n+1})) P_{n+1}$ with support in an isolated point $\theta_{n+1}^j$, which is of the form
\begin{equation}
\label{e:phi5}
  a(\theta_n^i), \ {\rm where} \ \theta_n^i \ \text{is an isolated point in} \ Z_n;
\end{equation}
\item an element of $P_{n+1} M_{\infty}(C(Z_{n+1})) P_{n+1}$, which is of the form
\begin{equation}
\label{e:phi6}
  (a_{ij} \cdot q)_{ij}, \ {\rm where} \ q \ \text{is a line bundle over} \ Z_{n+1}, \ {\rm and} \ (a_{ij}) = a(\theta_n^i);
\end{equation}
\item an element of $P_{n+1} M_{\infty}(C(Z_{n+1})) P_{n+1}$ of the form
\begin{equation}
\label{e:phi7}
  a \circ \lambda,
\end{equation}
where $\lambda: \: Z_{n+1} \to Z_n$ is a continuous map whose image is only contained in one wedge summand of $Z_n$ (see our constructions in \S~\ref{s:AH}).
\end{itemize}
\setlength{\parindent}{0cm} \setlength{\parskip}{0.5cm}

Note that in \eqref{e:phi2} and \eqref{e:phi4}, we identify $P_n(\theta_n^{\bm{i}}) M_{\infty} P_n(\theta_n^{\bm{i}})$ with $M_{[n,\bm{i}]}$ via a fixed isomorphism.

Let $P^{\bm{a}} \in M(A_{n+1})$ be projections, with $\sum_{\bm{a}} P^{\bm{a}} = 1$, giving rise to the diagonal form in \eqref{e:FormOfConnMaps}, and let $\varphi^{\bm{a}}$ be the homomorphism $A_n \to P_{\bm{a}} A_{n+1} P_{\bm{a}}, \, x \ma P^{\bm{a}} u \varphi(a) u^* P^{\bm{a}}$. Since each of the $P^{\bm{a}}$ either lies in $C([0,1],E_{n+1}^q)$ or $F_{n+1}$, we have $\img(\varphi^{\bm{a}}) \subseteq P^{\bm{a}} C([0,1],E_{n+1}^q) P^{\bm{a}}$ or $\img(\varphi^{\bm{a}}) \subseteq P^{\bm{a}} F_{n+1} P^{\bm{a}}$. Then both maps in \eqref{e:ConnMapPartOne}, \eqref{e:ConnMapPartTwo} are of the form $u^* (\bigoplus_{\bm{a}} \varphi^{\bm{a}}) u$. The unitary $u$ is a permutation matrix for the map in \eqref{e:ConnMapPartOne} and is a unitary in $C([0,1],E_{n+1})$ such that $u(0)$ and $u(1)$ are permutation matrices for the map in \eqref{e:ConnMapPartTwo}.

\bremark
\label{R:unitaries}
Let us write $C_n \defeq C([0,1],E_n)$ and $u_{n+1} \in C_{n+1}$ for the unitary for the map in \eqref{e:ConnMapPartTwo}. The only reason we need $u_{n+1}$ is to ensure that we send $(f,a)$ to an element satisfying the right boundary conditions at $t=0$ and $t=1$. For this, only the values $u_{n+1,t} \defeq u_{n+1}(t)$ at $t \in \gekl{0,1}$ matter. Therefore, by an iterative process, we can change $\beta_t$ in order to arrange $u_{n+1}=1$ for the map in \eqref{e:ConnMapPartTwo}: First of all, it is easy to see that $\varphi_n$ extends uniquely to a homomorphism $\Phi_n: \: C_n \oplus F_n \to C_{n+1} \oplus F_{n+1}$. Let us write $\Phi_n^C$ and $\Phi_n^F$ for the composites
$$
  C_n \oplus F_n \overset{\Phi_n}{\lori} C_{n+1} \oplus F_{n+1} \onto C_{n+1} \ {\rm and} \ C_n \oplus F_n \overset{\Phi_n}{\lori} C_{n+1} \oplus F_{n+1} \onto F_{n+1}.
$$
As $\varphi_n$ sends strictly positive elements to strictly positive elements, $\Phi_n$ is unital. Now, for all $n$, let $\Lambda_n(t) \subseteq [0,1]$ be a finite set such that for all $(f_n,a_n) \in A_n$ with $\varphi(f_n,a_n) = (f_{n+1},a_{n+1}) \in A_{n+1}$, $f_n \vert_{\Lambda_n(t)} \equiv 0$ implies $f_{n+1}(t) = 0$. In other words, $\Lambda_n(t)$ are the evaluation points for $f_{n+1}(t)$. Similarly, let $T_n \subseteq (0,1)$ be such that for all $(f_n,a_n) \in A_n$ with $\varphi(f_n,a_n) = (f_{n+1},a_{n+1}) \in A_{n+1}$, $f_n \vert_{T_n} \equiv 0$ and $a_n = 0$ imply $a_{n+1} = 0$. Now we choose inductively on $n$ unitaries $v_n \in U(C_n)$ and $u_{n+1} \in U(C_{n+1})$ such that, for all $n$, $v_n(s) = 1$ for all $s \in (\Lambda_n(0) \cup \Lambda_n(1) \cup T_n) \setminus \gekl{0,1}$, $u_{n+1}(t) = u_{n+1,t}$ for $t \in \gekl{0,1}$, and $v_{n+1} = \Phi_n^C(v_n,1) u_{n+1}^*$: Simply start with $v_1 \defeq 1$, and if $v_n$ and $u_n$ have been chosen, choose $u_{n+1} \in U(C_{n+1})$ such that $u_{n+1}(t) = u_{n+1,t}$ for all $t \in \gekl{0,1}$ and $u_{n+1}(s) = \Phi_n^C(v_n,1)(s)$ for all $s \in (\Lambda_n(0) \cup \Lambda_n(1) \cup T_n) \setminus \gekl{0,1}$, and set $v_{n+1} \defeq \Phi_n^C(v_n,1) u_{n+1}^*$. If we now take this $u_{n+1}$ for the map in \eqref{e:ConnMapPartTwo} giving rise to $\varphi_n$ and $\Phi_n$, then we obtain a commutative diagram
$$
  \xymatrix{
  C_n \oplus F_n \ar[d]^{\cong}_{(v_n^* \sqcup v_n) \oplus \id}  \ar[rr]^{\Phi_n} & & C_{n+1} \oplus F_{n+1} \ar[d]_{\cong}^{(v_{n+1}^* \sqcup v_{n+1}) \oplus \id}
  \\
  C_n \oplus F_n \ar[rr]^{(u_{n+1},1) \Phi_n (u_{n+1},1)^*} & \ \ \ \ \ \ \ \ & C_{n+1} \oplus F_{n+1}
  }
$$
which restricts to
$$
  \xymatrix{
  A_n \ar[r]^{\Phi_n} \ar[d]^{\cong}_{(v_n^* \sqcup v_n) \oplus \id} & A_{n+1} \ar[d]_{\cong}^{(v_{n+1}^* \sqcup v_{n+1}) \oplus \id}
  \\
  \bar{A}_n \ar[r]^{\bar{\varphi}_n} & \bar{A}_{n+1}
  }  
$$
where the unitary $\bar{u}_{n+1}$ for the map in \eqref{e:ConnMapPartTwo} for $\bar{\varphi}_n$ is now trivial, $\bar{u}_{n+1} = 1$, and $\bar{A}_n$ is of the same form \eqref{e:A} as $A_n$, with $\bar{\beta}_t = v_n(t)^* \beta_t v_n(t)$ of the same form as $\beta_t$ for $t = 0, 1$ (the point being that $v_n(t)$ is a permutation matrix). Obviously, we have $\ilim \gekl{\bar{A}_n;\bar{\varphi}_n} \cong \ilim \gekl{A_n;\varphi_n}$.
\eremark

\bremark
Note that the construction described in \S~\ref{ss:SummConst} also encompasses (a slight modification of) the C*-algebra construction in \cite[\S~6]{GLII}. (In particular, one obtains model algebras of rational generalized tracial rank one, in the sense of \cite{GLII}.)
\eremark

\section{Inductive limits and Cartan pairs revisited}
\label{s:LimCartan}

In this section, we improve the main result in \cite[\S~3]{BL17} and give a C*-algebraic interpretation. Let us first recall \cite[Theorem~3.6]{BL17}. We use the same notations and definitions as in \cite{Ren,BL17}. We start with the following
\bremark
We can drop the assumptions of second countability for groupoids and separability for C*-algebras in \cite{Ren} if we replace \an{topologically principal} by \an{effective} throughout. In other words, given a twisted \'{e}tale effective groupoid $(G,\Sigma)$, i.e., a twisted \'{e}tale groupoid $(G,\Sigma)$ where $G$ is effective (not necessarily second countable), $(C^*_r(G,\Sigma),C_0(G^{(0)}))$ is a Cartan pair; and conversely, given a Cartan pair $(A,B)$ (where $A$ is not necessarily separable), the Weyl twist $(G(A,B),\Sigma(A,B))$ from \cite{Ren} is a twisted \'{e}tale effective groupoid. These constructions are inverse to each other, i.e., there are canonical isomorphisms $(G,\Sigma) \cong (G(C^*_r(G,\Sigma),C_0(G^{(0)})),\Sigma(C^*_r(G,\Sigma),C_0(G^{(0)})))$ (provided by \cite[4.13, 4.15, 4.16]{Ren}) and $(A,B) \cong (C^*_r(G(A,B),\Sigma(A,B)),C_0(G(A,B)^{(0)}))$ (provided by \cite[5.3, 5.8, 5.9]{Ren}). Similarly, everything in \cite[\S~3]{BL17} works without the assumption of second countability. In particular, \cite[Theorem~3.6]{BL17} holds for general twisted \'{e}tale groupoids if we replace \an{topologically principal} by \an{effective}. This is why in this section, we formulate everything for twisted \'{e}tale effective groupoids and general Cartan pairs. In our applications later on, however, we will only consider second countable groupoids and separable C*-algebras.
\eremark

Now suppose that $(A_n,B_n)$ are Cartan pairs, let $(G_n,\Sigma_n)$ be their Weyl twists, and set $X_n \defeq G_n^{(0)}$. Let $\varphi_n: \: A_n \to A_{n+1}$ be injective *-homomorphisms. Assume that there are twisted groupoids $(H_n,T_n)$, with $Y_n \defeq H_n^{(0)}$, together with twisted groupoid homomorphisms $(i_n,\imath_n): \: (H_n,T_n) \to (G_{n+1},T_{n+1})$ and $(\dot{p}_n,p_n): \: (H_n,T_n) \to (G_n,T_n)$ such that $i_n: \: H_n \to G_{n+1}$ is an embedding with open image, and $\dot{p}_n: \: H_n \to G_n$ is surjective, proper, and fibrewise bijective (i.e., for every $y \in Y_n$, $\dot{p}_n \vert_{(H_n)_y}$ is a bijection onto $(G_n)_{\dot{p}_n(y)}$). Suppose that $\varphi_n = (\imath_n)_* \circ p_n^*$ for all $n$. Further assume that condition (LT) is satisfied, i.e., for every continuous section $\rho: \: U \to \rho(U)$ for the canonical projection $\Sigma_n \onto G_n$, where $U \subseteq G_n$ is open, there is a continuous section $\ti{\rho}: \: \dot{p}_n^{-1}(U) \to \ti{\rho}(\dot{p}_n^{-1}(U))$ for the canonical projection $T_n \onto H_n$ such that $\ti{\rho}(\dot{p}_n^{-1}(U)) \subseteq \dot{p}_n^{-1}(\rho(U))$ and $p_n \circ \ti{\rho} = \rho \circ \dot{p}_n$. Also assume that condition (E) is satisfied, i.e., for every continuous section $t: \: U \to t(U)$ for the source map of $G_n$, where $U \subseteq X_n$ and $t(U) \subseteq G_n$ are open, there is a continuous section $\ti{t}: \: \dot{p}_n^{-1}(U) \to \ti{t}(\dot{p}_n^{-1}(U))$ for the source map of $H_n$ such that $\ti{t}(\dot{p}_n^{-1}(U)) \subseteq \dot{p}_n^{-1}(t(U))$ and $\dot{p}_n \circ \ti{t} = t \circ \dot{p}_n$.

In this situation, define 
\begin{eqnarray}
\label{e:ilimGSigma}
  && \Sigma_{n,0} \defeq \Sigma_n \ {\rm and} \ \Sigma_{n,m+1} \defeq p_{n+m}^{-1}(\Sigma_{n,m}) \subseteq T_{n+m} \ {\rm for \ all} \ n \ {\rm and} \ m = 0, 1, \dotsc,\\
  && G_{n,0} \defeq G_n \ {\rm and} \ G_{n,m+1} \defeq \dot{p}_{n+m}^{-1}(G_{n,m}) \subseteq H_{n+m} \ {\rm for \ all} \ n \ {\rm and} \ m = 0, 1, \dotsc, \nonumber \\
  && \bar{\Sigma}_n \defeq \plim_m \gekl{\Sigma_{n,m}; p_{n+m}} \ {\rm and} \ \bar{G}_n \defeq \plim_m \gekl{G_{n,m}; \dot{p}_{n+m}} \ {\rm for \ all} \ n. \nonumber
\end{eqnarray}
Then \cite[Theorem~3.6]{BL17} tells us that
\begin{enumerate}
\item[(a)] $(\bar{G}_n,\bar{\Sigma}_n)$ are twisted groupoids, and $(i_n,\imath_n)$ induce twisted groupoid homomorphisms $(\bar{i}_n,\bar{\imath}_n): \: (\bar{G}_n,\bar{\Sigma}_n) \to (\bar{G}_{n+1},\bar{\Sigma}_{n+1})$ such that $\bar{i}_n$ is an embedding with open image for all $n$, and 
\begin{equation}
\label{e:barSigma,G}
  \bar{\Sigma} \defeq \ilim \gekl{\bar{\Sigma}_n; \bar{\imath}_n} \ {\rm and} \ \bar{G} \defeq \ilim \gekl{\bar{G}_n; \bar{i}_n}
\end{equation}
defines a twisted \'{e}tale groupoid $(\bar{G},\bar{\Sigma})$, 
\item[(b) \& (c)] $(\ilim \gekl{A_n;\varphi_n},\ilim \gekl{B_n;\varphi_n})$ is a Cartan pair whose Weyl twist is given by $(\bar{G},\bar{\Sigma})$.
\end{enumerate}
\bremark
\label{R:LimCDiag}
It is clear that the proof of \cite[Theorem~3.6]{BL17} shows that if all $B_n$ are C*-diagonals, i.e., all $G_n$ are principal, then $\bar{G}$ is principal, i.e., $\ilim \gekl{B_n;\varphi_n}$ is a C*-diagonal.
\eremark

It turns out that conditions (LT) and (E) are redundant.
\blemma
\label{L:LT,E}
In the situation above, conditions (LT) and (E) are automatically satisfied.
\elemma
\setlength{\parindent}{0cm} \setlength{\parskip}{0cm}

\bproof
To prove condition (LT), let $\rho: \: U \to \rho(U)$ be a continuous section for the canonical projection $\pi_n: \: \Sigma_n \onto G_n$, where $U \subseteq G_n$ is open. Let $\pi_{n+1}: \: \Sigma_{n+1} \onto G_{n+1}$ be the canonical projection. Then $\pi_{n+1} \vert_{p_n^{-1}(\rho(U))}: \: p_n^{-1}(\rho(U)) \to \dot{p}_n^{-1}(U)$ is bijective. Indeed, given $\tau_1, \tau_2 \in p_n^{-1}(\rho(U))$ with $\pi_{n+1}(\tau_1) = \pi_{n+1}(\tau_2) \eqdef \eta \in H_n$, we must have $\tau_2 = z \cdot \tau_1$ for some $z \in \Tz$. Also, $\pi_n(p_n(\tau_1)) = \dot{p}_n(\eta) = \pi_n(p_n(\tau_1))$. As $\pi_n \vert_{\rho(U)}: \: \rho(U) \to U$ is bijective (with inverse $\rho$), we deduce $p_n(\tau_1) = p_n(\tau_2)$. Hence $p_n(\tau_1) = p_n(\tau_2) = z \cdot p_n(\tau_1)$, which implies $z = 1$, i.e., $\tau_2 = \tau_1$. This proves injectivity, and surjectivity is easy to see. As $\pi_{n+1}$ is open, $\ti{\rho} \defeq (\pi_{n+1} \vert_{p_n^{-1}(\rho(U))})^{-1}: \: \dot{p}_n^{-1}(U) \to p_n^{-1}(\rho(U))$ is the continuous section we are looking for.

To verify (E), let $t: \: U \to t(U)$ be a continuous section for the source map $s_n$ of $G_n$, where $U \subseteq X_n$ and $t(U) \subseteq G_n$ are open. Let $s_{n+1}$ be the source map of $H_n$. Then $s_{n+1} \vert_{\dot{p}_n^{-1}(t(U))}: \: \dot{p}_n^{-1}(t(U)) \to \dot{p}_n^{-1}(U)$ is bijective. Indeed, given $\eta_1, \eta_2 \in \dot{p}_n^{-1}(t(U))$ with $s_{n+1}(\eta_1) = s_{n+1}(\eta_2) \eqdef y \in Y_n$, we must have $s_n(\dot{p}_n(\eta_1)) = \dot{p}_n(y) = s_n(\dot{p}_n(\eta_2))$. As $s_n \vert_{t(U)}: \: t(U) \to U$ is bijective (with inverse $t$), we deduce $\dot{p}_n(\eta_1) = \dot{p}_n(\eta_2)$. Since $\dot{p}_n$ is fibrewise bijective, this implies $\eta_1 = \eta_2$. This proves injectivity, and surjectivity is easy to see. As $\dot{p}_n^{-1}(t(U))$ is open and $s_{n+1}$ is open, $\ti{t} \defeq (s_{n+1} \vert_{\dot{p}_n^{-1}(t(U))})^{-1}: \: \dot{p}_n^{-1}(U) \to \dot{p}_n^{-1}(t(U))$ is the continuous section we are looking for.
\eproof
\setlength{\parindent}{0cm} \setlength{\parskip}{0.5cm}

Let us now determine which *-homomorphisms are of the form $\imath_* \circ p^*$. Let $(A,B)$ and $(\hat{A},\hat{B})$ be Cartan pairs with normalizers $N \defeq N_A(B)$, $\hat{N} \defeq N_{\hat{A}}(\hat{B})$ and faithful conditional expectations $P: \: A \onto B$, $\hat{P}: \: \hat{A} \onto \hat{B}$. Let $(G,\Sigma)$ and $(\hat{G},\hat{\Sigma})$ be the Weyl twists of $(A,B)$ and $(\hat{A},\hat{B})$. Suppose that $\varphi: \: A \to \hat{A}$ is an injective *-homomorphism.
\bprop
\label{P:phi=ip}
The following are equivalent:
\setlength{\parindent}{0cm} \setlength{\parskip}{0cm}

\begin{enumerate}
\item[(i)] $\varphi(B) \subseteq \hat{B}$, $\varphi(N) \subseteq \hat{N}$, $\hat{P} \circ \varphi = \varphi \circ P$;
\item[(ii)] There exists a twisted \'{e}tale effective groupoid $(H,T)$ and twisted groupoid homomorphisms $(i,\imath): \: (H,T) \to (\hat{G},\hat{\Sigma})$, $(\dot{p},p): \: (H,T) \to (G,\Sigma)$, where $i$ is an embedding with open image and $\dot{p}$ is surjective, proper and fibrewise bijective, such that $\varphi = \imath_* \circ p^*$.
\end{enumerate}
\eprop
\setlength{\parindent}{0cm} \setlength{\parskip}{0cm}

\bproof
(ii) $\Rarr$ (i): It is easy to see that $(\imath_* \circ p^*)(B) \subseteq \hat{B}$. Given an open bisection $S$ of $G$, $\dot{p}^{-1}(S)$ is an open bisection of $H$, and then $i(\dot{p}^{-1}(S))$ is an open bisection of $\hat{G}$. Therefore, $(\imath_* \circ p^*)(N) \subseteq \hat{N}$. Finally, we have $\hat{P} \circ (\imath_* \circ p^*) = (\imath_* \circ p^*) \circ P$ because $\dot{p}^{-1}(G^{(0)}) = H^{(0)}$.
\setlength{\parindent}{0cm} \setlength{\parskip}{0.5cm}

(i) $\Rarr$ (ii): Let $\breve{B}$ be the ideal of $\hat{B}$ generated by $\varphi(B)$, and $\breve{A} \defeq C^*(\varphi(A),\breve{B})$. Then $(\breve{A},\breve{B})$ is a Cartan pair: It is clear that $\breve{B}$ contains an approximate unit for $\breve{A}$. To see that $\breve{B}$ is maximal abelian, take $a \in \breve{A} \cap (\breve{B})'$. Let $b \in \hat{B}$, and take an approximate unit $(h_{\lambda}) \subseteq \breve{B}$ for $\breve{A}$. Then $ba = \lim_{\lambda} b h_{\lambda} a = \lim_{\lambda} a b h_{\lambda} = \lim_{\lambda} a h_{\lambda} b = ab$. Hence $a \in \breve{A} \cap (\hat{B})' = \breve{A} \cap \hat{B} = \breve{B}$ (the last equality holds because $\breve{B}$ contains an approximate unit for $\breve{A}$, and $\breve{B} \cdot \hat{B} \subseteq \breve{B}$). This shows $\breve{A} \cap (\breve{B})' = \breve{B}$. Moreover, we have $\varphi(N) \subseteq \breve{N} \defeq N_{\breve{A}}(\breve{B})$: Let $n \in \varphi(N)$, $b \in \breve{B}$, and $(h_{\lambda}) \subseteq B$ be an approximate unit for $A$. Then $nbn^* \in \hat{B}$ as $n \in \varphi(N) \subseteq \hat{N}$, and thus $nbn^* = \lim_{\lambda} \varphi(h_{\lambda}) n b n^* \subseteq \overline{\varphi(B) \cdot \hat{B}} \subseteq \breve{B}$. Finally, it is clear that $\breve{P} \defeq \hat{P} \vert_{\breve{A}}$ is a faithful conditional expectation onto $\breve{B}$.

Let $(H,T)$ be the Weyl twist attached to $(\breve{A},\breve{B})$, and write $X \defeq G^{(0)}$, $Y \defeq H^{(0)}$ and $\hat{X} \defeq \hat{G}^{(0)}$. It is easy to see that $\breve{N} \subseteq \hat{N}$. Hence we may define maps
$$
  i: \: H \to \hat{G}, \, [x,\alpha_n,y] \ma [x,\alpha_n,y] \ \ \ {\rm and} \ \ \ \imath: \: T \to \hat{\Sigma}, \, [x,n,y] \ma [x,n,y], \ {\rm for} \ n \in \breve{N}. 
$$
Clearly, $i$ and $\imath$ are continuous groupoid homomorphisms. $i$ is injective since $[x,\alpha_n,y] = [x',\alpha_{n'},y']$ in $\hat{G}$ implies $x = x'$, $y = y'$ and $\alpha_n = \alpha_{n'}$ on a neighbourhood $U \subseteq \hat{X}$ of $y$, so that $\alpha_n = \alpha_{n'}$ on $U \cap Y$, which is a neighbourhood of $y$ in $Y$, and hence $[x,\alpha_n,y] = [x',\alpha_{n'},y']$ in $H$. The image of $i$ is given by $\bigcup_{n \in \breve{N}} \menge{[\alpha_n(y),\alpha_n,y]}{y \in \dom(n)}$ which is clearly open in $\hat{G}$. Finally, it is easy to see that we have a commutative diagram
$$
  \xymatrix{
  \breve{A} \ar[r]_{\subseteq} \ar[d]_{\cong} & \hat{A} \ar[d]^{\cong}  \\
  C^*_r(H,T) \ar[r]^{\imath_*} & C^*_r(\hat{G},\hat{\Sigma}),
  }  
$$
where the upper horizontal map is given by inclusion, and the vertical isomorphisms are as in \cite[Definition~5.4]{Ren}.

We now proceed to construct $(\dot{p},p)$. Since $A = C^*(N)$ and $\varphi(N) \subseteq \breve{N}$, it is easy to see that $\breve{A} = \clspan(\varphi(N) \cdot \breve{B})$. It follows that for every $\breve{n} \in \breve{N}$ and $y \in \dom(\breve{n})$, there is $n \in \varphi(N)$ such that $y \in \dom(n)$ and $[x,\breve{n},y] = [x,n,y]$ in $T$. Indeed, for $a \in \lspan(\varphi(N) \cdot \breve{B})$ it is clear that $a \equiv 0$ on $T \setminus \big( \bigcup_{n \in \varphi(N)} \menge{[\alpha_n(y),n,y]}{y \in \dom(n)} \big)$. As the latter set is closed in $T$, we must have $a \equiv 0$ on $T \setminus \big( \bigcup_{n \in \varphi(N)} \menge{[\alpha_n(y),n,y]}{y \in \dom(n)} \big)$ for all $a \in \breve{A}$. Hence $T = \bigcup_{n \in \varphi(N)} \menge{[\alpha_n(y),n,y]}{y \in \dom(n)}$. This observation allows us to define the maps
$$
  \dot{p}: \: H \to G, \, [x,\alpha_{\varphi(n)},y] \ma [\varphi^*(x),\alpha_n,\varphi^*(y)] \ \ \ {\rm and} \ \ \ p: \: T \to \Sigma, \, [x,\varphi(n),y) \ma [\varphi^*(x),n,\varphi^*(y)], \ {\rm for} \ n \in N,
$$
where $\varphi^*: \: Y \to X$ is the map dual to $B \to \breve{B}, \, b \ma \varphi(b)$ determined by $\varphi(b) = b \circ \varphi^*$ for all $b \in B$. Note that $\varphi^*$ exists since $\varphi(B)$ is full in $\breve{B}$. $p$ is well-defined because $[x,\varphi(m),y] = [x,\varphi(n),y]$ implies $\hat{P}(\varphi(n)^*\varphi(m))(y) > 0$, so that $P(n^*m)(\varphi^*(y)) = \varphi(P(n^*m))(y) = \hat{P}(\varphi(n)^* \varphi(m))(y) > 0$, which in turn yields $[\varphi^*(x), \alpha_m, \varphi^*(y)] = [\varphi^*(x), \alpha_n, \varphi^*(y)]$. Similarly, $\dot{p}$ is well-defined. Clearly, $(\dot{p},p)$ is a twisted groupoid homomorphism. As $\varphi$ is injective, $\varphi^*$ is surjective, so that $\dot{p}$ is surjective.
\setlength{\parindent}{0.5cm} \setlength{\parskip}{0cm}

To see that $\dot{p}$ is proper, let $K \subseteq G$ be compact. Given $n \in N$, write $U(n) \defeq \menge{[\alpha_n(y),\alpha_n,y]}{y \in \dom(n)}$ and $K(n) \defeq (s \vert_{U(n)})^{-1}(s(K))$. As $K$ is compact, there exists a finite set $\gekl{n_i} \subseteq N$ such that $K \subseteq \bigcup_i U(n_i)$, so that $K = \bigcup_i U(n_i) \cap K \subseteq \bigcup_i K(n_i)$. Now given $m \in N$, $\dot{p}([x,\alpha_{\varphi(m)},y]) \in K(n)$ implies $\varphi^*(y) \in s(K)$, i.e., $y \in (\varphi^*)^{-1}(s(K))$, $\dot{p}([x,\alpha_{\varphi(m)},y]) = [\varphi^*(x), \alpha_m, \varphi^*(y)] = [\varphi^*(x), \alpha_n, \varphi^*(y)]$, so that $P(n^*m)(\varphi^*(y)) \neq 0$, which yields $\hat{P}(\varphi(n)^* \varphi(m))(y) = \varphi(P(n^* m))(y) \neq 0$, thus $[x,\alpha_{\varphi(m)},y] = [x,\alpha_{\varphi(n)},y]$. Hence $\dot{p}^{-1}(K(n)) \subseteq \menge{[\alpha_{\varphi(n)}(y), \varphi(n), y]}{y \in (\varphi^*)^{-1}(s(K))} = (s \vert_{U(\varphi(n))})^{-1}((\varphi^*)^{-1}(s(K)) \eqdef \breve{K}(n)$. As $\varphi^*$ is proper, $\breve{K}(n)$ is compact for all $n \in N$. Hence $\dot{p}^{-1}(K) \subseteq \bigcup_i \dot{p}^{-1}(K(n_i)) \subseteq \bigcup_i \breve{K}(n_i)$ is a closed subset of a compact set, thus compact itself. 

Moreover, given $y \in Y$, $\dot{p}([w,\alpha_{\varphi(m)},y]) = \dot{p}([x,\alpha_{\varphi(n)},y])$ implies $[\varphi^*(w), \alpha_m, \varphi^*(y)] = [\varphi^*(x), \alpha_n, \varphi^*(y)]$, so that $\hat{P}(\varphi(n)^* \varphi(m))(y) = P(n^* m)(\varphi^*(y)) \neq 0$, so that $[w,\alpha_{\varphi(m)},y] = [x,\alpha_{\varphi(n)},y]$. This shows injectivity of $\dot{p} \vert_{H_y}$, and it is clear that $\dot{p}(H_y) = G_{\dot{p}(y)}$. Thus $\dot{p}$ is fibrewise bijective. 

Finally, it is easy to see that we have a commutative diagram
$$
  \xymatrix{
  A \ar[r]_{\varphi} \ar[d]_{\cong} & \breve{A} \ar[d]^{\cong}  \\
  C^*_r(G,\Sigma) \ar[r]^{p^*} & C^*_r(H,T),
  }  
$$
where the vertical isomorphisms are as in \cite[Definition~5.4]{Ren}.
\eproof
\setlength{\parindent}{0cm} \setlength{\parskip}{0.5cm}

\bremark
\label{R:full}
In Proposition~\ref{P:phi=ip}, $\varphi$ sends full elements to full elements if and only if we have $i(H^{(0)}) = \hat{G}^{(0)}$.
\eremark

Theorem~\ref{THM:CartanLim} now follows from \cite[Theorem~3.6]{BL17}, Lemma~\ref{L:LT,E}, Proposition~\ref{P:phi=ip} and Remark~\ref{R:LimCDiag}.

\bremark
\label{R:Weyl,full,diag}
The Weyl twist of $(\ilim \gekl{A_n; \varphi_n}, \ilim \gekl{B_n; \varphi_n})$ in the situation of Theorem~\ref{THM:CartanLim} is given by $(\bar{G},\bar{\Sigma})$ as given in \eqref{e:ilimGSigma} and \eqref{e:barSigma,G}.
\setlength{\parindent}{0.5cm} \setlength{\parskip}{0cm}

If, in Theorem~\ref{THM:CartanLim}, all $\varphi_n$ send full elements to full elements, then $G_{n,m+1}^{(0)} = H_{n+m}^{(0)} = G_{n+m+i}^{(0)}$ (where we identify $H_{n+m}^{(0)}$ with $i_{n+m}(H_{n+m}^{(0)})$), so that $\bar{G}_n^{(0)} = \plim_m \{ G_{n,m}^{(0)};\dot{p}_{n+m} \} \cong \plim \{ G_l^{(0)}; \dot{p}_l \}$ for all $n$, and thus $\bar{i}_n(\bar{G}_n^{(0)}) = \bar{G}_{n+1}^{(0)}$ for all $n$, which implies $\bar{G}^{(0)} \cong \plim \{ G_n^{(0)}; \dot{p}_n \}$.

If all $B_n$ in Theorem~\ref{THM:CartanLim} are C*-diagonals, i.e., all $G_n$ are principal, then $\bar{G}$ is principal.
\eremark

\section{Groupoid models}
\label{s:GPDModels}

\subsection{The building blocks}
\label{ss:BB}

We first present groupoid models for the building blocks that give rise to our AH-algebras (see \S~\ref{s:AH}). Let $Z$ be a second countable compact Hausdorff space and let $p_i \in M_{\infty}(C(Z))$ be a finite collection of line bundles over $Z$. Let $P = \sum_i{}^{\oplus} \ p_i \in M_{\infty}(C(Z))$. The following is easy to check:
\blemma
\label{L:F,F}
$\bigoplus_i p_i M_{\infty}(C(Z)) p_i$ is a Cartan subalgebra of $P M_{\infty}(C(Z)) P$.
\elemma
\setlength{\parindent}{0cm} \setlength{\parskip}{0cm}

Thus, by \cite[Theorem~5.9]{Ren}, there exists a twisted groupoid $(\dot{\cF},\cF)$ (the Weyl twist) such that 
$$
  (C^*_r(\dot{\cF},\cF), \, C_0(\dot{\cF}^{(0)})) \cong (P M_{\infty}(C(Z)) P, \, \bigoplus_i p_i M_{\infty}(C(Z)) p_i).
$$

Let us now describe $(\dot{\cF},\cF)$ explicitly. Let $R$ be the full equivalence relation on the finite set $\gekl{p_i}$ (just a set with the same number of elements as the number of line bundles). Let $\dot{\cF} = Z \times R$, which is a groupoid in the canonical way. For every $p_i$, let $T_i$ be a circle bundle over $Z$ such that $p_i = \Cz \times_{\Tz} T_i$. We form the circle bundles $T_j \cdot T_i^*$, which are given as follows: For each index $i$, let $\gekl{V_{i,a}}_a$ be an open cover of $Z$, and let $v_{i,a}$ be a trivialization of $T_i \vert_{V_{i,a}}$. We view $v_{i,a}$ as a continuous map $v_{i,a}: \: V_{i,a} \to M_{\infty}$ with values in partial isometries such that $v_{i,a}(z)$ has source projection $e_{11}$ and range projection $p_i(z)$, so that $v_{i,a}(z) = p_i(z) v_{i,a}(z) e_{11}$. Here $e_{11}$ is the rank one projection in $M_{\infty}$ which has zero entry everywhere except in the upper left $(1,1)$-entry, where the value is $1$. Then
$$
  T_j \cdot T_i^* = \Big( \coprod_{c,a} \Tz \times (V_{j,c} \cap V_{i,a}) \Big) \Big/ { }_{\sim}
$$
where we define $(z,x) \sim (z',x')$ if $x = x'$, and if $x \in V_{j,c} \cap V_{i,a}$, $x' \in V_{j,d} \cap V_{i,b}$, then $z' = v_{i,b} v_{j,d}^* v_{j,c} v_{i,a}^* z$. 
\setlength{\parindent}{0cm} \setlength{\parskip}{0.5cm}

We set
$$
  \cF \defeq \coprod_{j,i} T_j \cdot T_i^*.
$$
Note that $T_i \cdot T_i^*$ is just the trivial circle bundle $\Tz \times Z$. We define a multiplication on $\cF$: For $([z,x],(j,i))$ and $([z',x'],(j',i'))$ in $\cF$, we can only multiply these elements if $x = x'$ and $i = j'$. In that case, write $h \defeq i'$ and assume that $x \in V_{j,c} \cap V_{i,b}$ and $x' = x \in V_{i,b} \cap V_{h,a}$. Then we define the product as
$$
  ([z,x],(j,i)) \cdot ([z',x'],(j',i')) = ([zz',x],(j,h)).
$$
Moreover, $\cF$ becomes a twist of $\dot{\cF}$ via the map
$$
  \cF \to \dot{\cF}, \, T_j \cdot T_i^* \ni \sigma \ma (\pi(\sigma),(j,i)), \ {\rm where} \ \pi: \: T_j \cdot T_i^* \to Z \ \text{is the canonical projection}.
$$
It is now straightforward to check (compare \cite{Ren}) that the twisted groupoid $(\dot{\cF},\cF)$ is precisely the Weyl twist of
$(P M_{\infty}(C(Z)) P, \, \bigoplus_i p_i M_{\infty}(C(Z)) p_i)$. More precisely, we have the following
\blemma
We have a $C(Z)$-linear isomorphism $C^*_r(\dot{\cF},\cF) \cong P M_{\infty}(C(Z)) P$ sending $\ti{f} \in C_c(\dot{\cF},\cF)$ with $\supp(\ti{f}) \subseteq \rukl{V_{j,c} \cap V_{i,a}} \times \gekl{(j,i)} \subseteq \dot{\cF}$ to $f v_{j,c} v_{i,a}^*$, where $f \in C(Z)$ is determined by $\ti{f}(([z,x],(j,i)) = \bar{z} f(x)$. Moreover, this $C(Z)$-linear isomorphism identifies $C(\dot{\cF}^{(0)})$ with $\bigoplus_i p_i M_{\infty}(C(Z)) p_i$.
\elemma

Let us now fix $n$, and apply the result above to the homogeneous C*-algebra $F \defeq F_n$ from \S~\ref{ss:SummConst} to obtain a twisted groupoid $(\dot{\cF},\cF)$ such that $C^*_r(\dot{\cF},\cF) \cong F$. More precisely, we apply our construction above to the summand of $F$ corresponding to the component of $\theta_n^{\bm{i}}$. Note that in the construction above, all our line bundles satisfy
\begin{equation}
\label{e:ptheta=e}
  p_i(\theta_n^{\bm{i}}) = e_{11}
\end{equation}
because of \eqref{e:pXN=_}. For the other summands, it is easy to construct a groupoid model, as these are just matrix algebras, so that we can just take the full equivalence relation on finite sets. 

Now our goal is to present a groupoid model for the building block $A \defeq A_n$ in \S~\ref{ss:SummConst}. Let $\cR$ be an equivalence relation (on a finite set) such that $C^*(\cR) \cong E \defeq E_n$. Write $\cR = \coprod_p \cR^p$ for subgroupoids $\cR^p$ such that the isomorphism $C^*(\cR) \cong E$ restricts to isomorphisms $C^*(\cR^p) \cong E^p \defeq E_n^p$. Set $\dot{\cC} \defeq [0,1] \times \cR$. Then $C^*_r(\dot{\cC})$ is canonically isomorphic to $C \defeq C([0,1],E)$. Consider the trivial twist $\cC \defeq \Tz \times \dot{\cC}$ of $\dot{\cC}$. Clearly, we have $C^*_r(\dot{\cC} \amalg \dot{\cF}, \cC \amalg \cF) \cong C \oplus F$. 

For $t=0,1$ and $\beta_t$ as in \S~\ref{ss:SummConst}, write
$$
  F \overset{\beta_t}{\lori} E \onto E^p
$$
as the composition
\begin{equation}
\label{e:FtoE}
  F \to \bigoplus_l M_{n_l} \otimes \Cz^{I_l^p} \into E^p,
\end{equation}
where each of the components $F \to M_{n_l} \otimes \Cz^{I_l^p}$ of the first map is given by
$$
  a \ma   
    \begin{pmatrix}
    a(\theta^l) & & \\
    & a(\theta^l) & \\
    & & \ddots
    \end{pmatrix},
$$
with $\# I_l^p$ copies of $a(\theta^l)$ on the diagonal, and the components $M_{n_l} \otimes \Cz^{I_l^p} \into E^p$ of the second map are given by
\begin{equation}
\label{e:2CompBeta}
    \bm{x}
    \ma
    u^*
    \begin{pmatrix}
    \bm{x} & & & \\
     & 0 & & \\
     & & \ddots & \\
     & & & 0
    \end{pmatrix}
  u,
\end{equation}
where $u$ is a permutation matrix.

Let $E_t^p$ be the image of $\bigoplus_l M_{n_l} \otimes \Cz^{I_l^p}$ in $E$, and set $E_t \defeq \bigoplus_p E_t^p \subseteq E$, for $t=0,1$. Let $\cR_t^p \subseteq \cR^p$ be subgroupoids such that the identification $C^*(\cR^p) \cong E^p$ restricts to $C^*(\cR_t^p) \cong E_t^p$. Write $\cR_t \defeq \coprod_p \cR_t^p$, so that $C^*(\cR) \cong E$ restricts to $C^*(\cR_t) \cong E_t$. Let $\sigma_t^p$ be the groupoid isomorphism $\coprod_l \cR_l \times I_l^p \cong \cR_t^p$, given by a bijection of the finite unit space, corresponding to conjugation by the unitary $u$ in \eqref{e:2CompBeta}. Let $V_{i,a}$ and $v_{i,a}$ be as above (introduced after Lemma~\ref{L:F,F}). We now define a map $\bm{b}_t: \: \Tz \times (\gekl{t} \times \cR_t) \to \Sigma$ as follows: Given an index $l$ and $(j,i) \in \cR_l$, choose indices $a$ and $c$ such that $\theta^l \in V_{j,c} \cap V_{i,a}$. Then define
\begin{equation}
\label{e:zji--beta}
  z_{j,i} \defeq v_{j,c}(\theta^l) v_{i,a}(\theta^l)^* \in \Tz.
\end{equation}
Here, we are using \eqref{e:ptheta=e}. If $\theta^l$ is not the distinguished point $\theta_n^{\bm{i}}$, then we set $z_{j,i} = 1$. For $z \in \Tz$ and $h \in I_l^p$, set
$$
  \bm{b}_t(z,t,\sigma_t^p((j,i),h)) \defeq [z_{j,i},\theta^l] \in T_j \cdot T_i^* \subseteq \Sigma, \ \text{where we view} \ (z_{j,i},\theta^l) \ \text{as an element in} \ \Tz \times (V_{j,c} \cap V_{i,a}).
$$
Define
$$
  \check{\Sigma} \defeq \menge{x \in \cC \amalg \cF}{x = (z,t,\gamma) \in \Tz \times [0,1] \times \cR \ \Rarr \ \gamma \in \cR_t \ {\rm for} \ t = 0,1} \ \ \ {\rm and} \ \Sigma \defeq \check{\Sigma} / { }_{\sim}
$$
where $\sim$ is the equivalence relation on $\check{\Sigma}$ generated by $(z,t,\gamma) \sim \bm{b}_t(z,t,\gamma)$ for all $z \in \Tz$, $t=0,1$ and $\gamma \in \cR_t$. $\check{\Sigma}$ and $\Sigma$ are principal $\Tz$-bundles belonging to twisted groupoids, and we denote the underlying groupoids by $\check{G}$ and $G$.

By construction, the canonical projection and inclusion $\Sigma \twoheadleftarrow \check{\Sigma} \into \cC \amalg \cF$ induce on the level of C*-algebras
$$
  \xymatrix{
  C^*_r(\dot{\cC} \amalg \dot{\cF}, \cC \amalg \cF) \ar[r]_{\cong} & C \oplus F
  \\
  C^*_r(\check{G},\check{\Sigma}) \ar[u]_{\subseteq} \ar[r]_(.2){\cong} & \check{A} \defeq \menge{(f,a) \in C \oplus F}{f(t) \in E_t \ {\rm for} \ t = 0,1} \ar[u]^{\subseteq}
  \\
  C^*_r(G,\Sigma) \ar@{^{(}->}[u] \ar[r]^(.2){\cong} & A = \menge{(f,a) \in C \oplus F}{f(t) = \beta_t(a) \ {\rm for} \ t = 0,1} \ar[u]^{\subseteq}
  }
$$
In particular, $(G,\Sigma)$ is the desired groupoid model for our building block.

In what follows, it will be necessary to keep track of the index $n$, so that we will consider, for all $n$, twisted groupoids $(\dot{\cC}_n \amalg \dot{\cF}_n, \cC_n \amalg \cF_n)$, $(\check{G}_n,\check{\Sigma}_n)$, $(G_n,\Sigma_n)$ describing the C*-algebras $C_n \oplus F_n$, $\check{A}_n$ and $A_n$ as explained above. Moreover, for all $n$, let $B_n \subseteq A_n$ be the subalgebra corresponding to $C_0(G_n^{(0)})$ under the isomorphism $C^*_r(G_n,\Sigma_n) \cong A_n$.

\subsection{The connecting maps}
\label{ss:ConnMaps}

Let us now describe the connecting maps $\varphi_n: \: A_n \to A_{n+1}$ in the groupoid picture above. Let $P_{n+1}^{\bm{a}}$, $\varphi_n^{\bm{a}}$ be as in \S~\ref{ss:SummConst}, so that $\varphi_n = \bigoplus_{\bm{a}} \varphi_n^{\bm{a}}$ and $\img(\varphi_n^{\bm{a}}) \subseteq P_{n+1}^{\bm{a}} A_{n+1} P_{n+1}^{\bm{a}}$. Also, let $\Phi_n: \: C_n \oplus F_n \to C_{n+1} \oplus F_{n+1}$ be the extension of $\varphi_n$ as in Remark~\ref{R:unitaries}. Set $\Phi_n^{\bm{a}}: \: C_n \oplus F_n \to P_{n+1}^{\bm{a}}(C_{n+1} \oplus F_{n+1}) P_{n+1}^{\bm{a}}, \, x \ma P_{n+1}^{\bm{a}} \varphi_n^{\bm{a}}(x) P_{n+1}^{\bm{a}}$. We obtain $\check{\varphi}_n: \: \check{A}_n \to \check{A}_{n+1}$ and $\check{\varphi}_n^{\bm{a}}: \: \check{A}_n \to P_{n+1}^{\bm{a}} \check{A}_{n+1} P_{n+1}^{\bm{a}}$ by restricting $\Phi_n$ and $\Phi_n^{\bm{a}}$. Set
\begin{eqnarray*}
  (C \oplus F)[\Phi_n] &\defeq& \Big\{ x \in C_{n+1} \oplus F_{n+1}: \: x = \sum_{\bm{a}} P_{n+1}^{\bm{a}} x P_{n+1}^{\bm{a}} \Big\},\\
  \check{A}[\check{\varphi}_n^{\bm{a}}] &\defeq& \img(\check{\varphi}_n^{\bm{a}}),\\
  \check{A}[\varphi_n] &\defeq& \Big\{ x \in \check{A}_{n+1}: \: x = \sum_{\bm{a}} P_{n+1}^{\bm{a}} x P_{n+1}^{\bm{a}}, \, P_{n+1}^{\bm{a}} x P_{n+1}^{\bm{a}} \in \check{A}[\check{\varphi}_n^{\bm{a}}] \Big\},\\
  A[\varphi_n] &\defeq& A_{n+1} \cap \check{A}[\check{\varphi}_n].
\end{eqnarray*}
Note that $\check{A}[\check{\varphi}_n^{\bm{a}}] = P_{n+1}^{\bm{a}} F_{n+1} P_{n+1}^{\bm{a}}$ if $P_{n+1}^{\bm{a}} \in F_{n+1}$ and $\check{A}[\check{\varphi}_n^{\bm{a}}] = \menge{x \in P_{n+1}^{\bm{a}} \check{A}_{n+1} P_{n+1}^{\bm{a}}}{x(t) \in \img(\ev_t \circ \check{\varphi}_n^{\bm{a}}) \ {\rm for} \ t = 0,1}$ if $P_{n+1}^{\bm{a}} \in C_{n+1}$.

Let $\cT_n$ be the open subgroupoid of $\cC_{n+1} \amalg \cF_{n+1}$, with $\dot{\cT}_n \subseteq \dot{\cC}_{n+1} \amalg \dot{\cF}_{n+1}$ correspondingly, such that $C^*_r(\dot{\cC}_{n+1} \amalg \dot{\cF}_{n+1}, \cC_{n+1} \amalg \cF_{n+1}) \cong C_{n+1} \oplus F_{n+1}$ restricts to $C^*_r(\dot{\cT}_n,\cT_n) \cong (C \oplus F)[\Phi_n]$. Similarly, let $\check{T}_n$ be the open subgroupoid of $\check{\Sigma}_{n+1}$, with $\check{H}_n \subseteq \check{G}_{n+1} $ correspondingly, such that $C^*_r(\check{G}_{n+1}, \check{\Sigma}_{n+1}) \cong \check{A}_{n+1}$ restricts to $C^*_r(\check{H}_n,\check{T}_n) \cong \check{A}[\check{\varphi}_n]$. For $\eta \in \check{T}_n$ and $\eta' \in \check{\Sigma}_{n+1}$, $\eta \sim \eta'$ implies that $\eta'$ lies in $\check{T}_n$. It follows that $T_n = \check{T}_n / { }_{\sim}$ is an open subgroupoid of $\Sigma_{n+1}$. Define $H_n = \check{H}_n / { }_{\sim}$ in a similar way. By construction, the commutative diagram at the groupoid level
$$
  \xymatrix{
  \cC_n \amalg \cF_n & \cT_n \ar[r]_(.3){\subseteq} & \cC_{n+1} \amalg \cF_{n+1}
  \\
  \check{\Sigma}_n \ar[u]_{\subseteq} \ar@{->>}[d] & \check{T}_n \ar[r]^{\subseteq} \ar[u]_{\subseteq} \ar@{->>}[d] & \check{\Sigma}_{n+1} \ar[u]^{\subseteq} \ar@{->>}[d]
  \\
  \Sigma_n & T_n \ar[r]^{\subseteq} & \Sigma_{n+1}
  }
$$
induces at the C*-level
$$
  \xymatrix{
  C_n \oplus F_n & (C \oplus F)[\Phi_n] \ar[r]_{\subseteq} & C_{n+1} \oplus F_{n+1}
  \\
  \check{A}_n \ar[u]_{\subseteq} & \check{A}[\check{\varphi}_n] \ar[u]_{\subseteq} \ar[r]_{\subseteq} & \check{A}_{n+1} \ar[u]^{\subseteq}
  \\
  A_n \ar[u]_{\subseteq} & A[\varphi_n] \ar[u]_{\subseteq} \ar[r]^{\subseteq} & A_{n+1} \ar[u]^{\subseteq}
  }
$$
Let $\check{T}_n = \coprod_{\bm{a}} \check{T}_n^{\bm{a}}$ and $\check{H}_n = \coprod_{\bm{a}} \check{H}_n^{\bm{a}}$ be the decompositions into subgroupoids such that the identification $C^*_r(\check{H}_n, \check{T}_n) \cong \check{A}[\check{\varphi}_n] \subseteq \check{A}_{n+1}$ restricts to $C^*_r(\check{H}_n^{\bm{a}}, \check{T}_n^{\bm{a}}) \cong \check{A}[\check{\varphi}_n^{\bm{a}}]$. For fixed $n$ and every $\varphi^{\bm{a}} = \varphi^{\bm{a}}_n$ from our list in \S~\ref{ss:SummConst}, we now construct a map $p^{\bm{a}}: \: \check{T}_n^{\bm{a}} \to \Sigma_n$ such that
$$
  \xymatrix{
  C^*_r(G_n, \Sigma_n) \ar[d]^{\cong} \ar[r]_{C^*(p^{\bm{a}})} & C^*_r(\check{H}_n^{\bm{a}}, \check{T}_n^{\bm{a}}) \ar[d]_{\cong}
  \\
  A_n \ar[r]^{\varphi^{\bm{a}}} & \check{A}[\check{\varphi}_n^{\bm{a}}]
  }
$$
commutes. 

Recall that $\check{\Sigma}_n \subseteq \cC_n \amalg \cF_n = (\Tz \times [0,1] \times \cR_n) \amalg \cF_n$. Also, we denote the canonical projection $\cF_n \onto \dot{\cF}_n$ by $\sigma \ma \dot{\sigma}$.

\begin{itemize}
\item For $\varphi^{\bm{a}}$ as in \eqref{e:phi1}, let $p^{\bm{a}}$ be the composite
\begin{align*}
  \check{T}_n^{\bm{a}} \cong \Tz \times [0,1] \times \cR_n &\to \check{\Sigma}_n \overset{q}{\lori} \Sigma_n,\\
  (z,t,\gamma) &\ma (z,\lambda(t),\gamma)
\end{align*} 
where we note that the first map has image in $\check{\Sigma}_n$, so that we an apply the quotient map $q: \: \check{\Sigma}_n \onto \Sigma_n$.
\item For $\varphi^{\bm{a}}$ as in \eqref{e:phi2}, let $p^{\bm{a}}$ be the composite
\begin{align}
\label{e:ztgamma}
  \check{T}_n^{\bm{a}} \cong \Tz \times [0,1] \times \cR_n &\to \cF_n \overset{q}{\lori} \Sigma_n,\\
  (z,t,\gamma) &\ma z \cdot \sigma(t,\gamma) \nonumber
\end{align}
where $\sigma$ is a continuous groupoid homomorphism such that $\dot{\sigma}(t,\gamma) = (x(t),\gamma)$. For $x(t) = \theta^l \in \gekl{\theta_n^i}$ and $\gamma = (j,i)$, write
\begin{equation}
\label{e:sigmatgamma}
  \sigma(t,\gamma) = [z_{j,i},\theta^l], 
\end{equation}
which has to match up with \eqref{e:zji--beta}.
\item For $\varphi^{\bm{a}}$ as in \eqref{e:phi3}, let $p^{\bm{a}}$ be the composite
\begin{align*}
  \check{T}_n^{\bm{a}} \cong \Tz \times \{ \theta_{n+1}^j \} \times \cR_n &\to \check{\Sigma}_n \overset{q}{\lori} \Sigma_n,\\
  (z,\theta_{n+1}^j,\gamma) &\ma (z,\bm{t},\gamma).
\end{align*}
\item For $\varphi^{\bm{a}}$ as in \eqref{e:phi4}, let $p^{\bm{a}}$ be the composite
\begin{align*}
  \check{T}_n^{\bm{a}} \cong \Tz \times \{ \theta_{n+1}^j \} \times (\dot{\cF}_n)^x_x &\to \cF_n \overset{q}{\lori} \Sigma_n,\\
  (z,\theta_{n+1}^j,\gamma) &\ma z \cdot \sigma(\gamma),
\end{align*}
where $\sigma: \: (\dot{\cF}_n)^x_x \to \cF_n$ is a groupoid homomorphism with $\dot{\sigma}(\gamma) = (x,\gamma)$ matching up with $\sigma$ in \eqref{e:ztgamma}.
\item For $\varphi^{\bm{a}}$ as in \eqref{e:phi5}, let $p^{\bm{a}}$ be the composite
\begin{align*}
  \check{T}_n^{\bm{a}} \cong \Tz \times \{ \theta_{n+1}^j \} \times (\dot{\cF}_n)^{\theta_n^i}_{\theta_n^i} &\to \cF_n \overset{q}{\lori} \Sigma_n,\\
  (z,\theta_{n+1}^j,\gamma) &\ma (z,\theta_n^i,\gamma).
\end{align*}
\item For $\varphi^{\bm{a}}$ as in \eqref{e:phi6}, let $p^{\bm{a}}$ be the composite
\begin{align*}
  \check{T}_n^{\bm{a}} \cong \Tz \times Z_{n+1} \times (\dot{\cF}_n)^{\theta_n^i}_{\theta_n^i} \onto \Tz \times (\dot{\cF}_n)^{\theta_n^i}_{\theta_n^i} &\to \cF_n \overset{q}{\lori} \Sigma_n,\\
  (z,\gamma) &\ma z \cdot \sigma(\gamma),
\end{align*}
where $\sigma: \: (\dot{\cF}_n)^{\theta_n^i}_{\theta_n^i} \to \cF_n$ is a groupoid homomorphism with $\dot{\sigma}(\gamma) = (\theta_n^i,\gamma)$ matching up with \eqref{e:zji--beta}, just as \eqref{e:sigmatgamma}.
\item For $\varphi^{\bm{a}}$ as in \eqref{e:phi7}, we have $C^*_r(\check{H}_n^{\bm{a}}, \check{T}_n^{\bm{a}}) \cong \rukl{\sum_i \lambda^*(p_i)} \cdot F_{n+1} \cdot \rukl{\sum_i \lambda^*(p_i)}$, where $p_i$ are the line bundles such that $P_n = \sum_i p_i$ (see \S~\ref{s:AH} and \S~\ref{ss:BB}), and $p^{\bm{a}}$ is the composite
\begin{align*}
  \check{T}_n^{\bm{a}} &\to \cF_n \overset{q}{\lori} \Sigma_n,\\
  [z,x] &\ma [z,\lambda(x)],
\end{align*}
with $(z,x) \in \Tz \times \lambda^{-1}(V_{i,a})$ and $(z,\lambda(x)) \in \Tz \times V_{i,a}$, where for a given open cover $V_{i,a}$ and trivialization $v_{i,a}$ for $\cF_n$, we choose the open cover $\lambda^{-1}(V_{i,a})$ and trivialization $v_{i,a} \circ \lambda$ for $\check{T}_n^{\bm{a}}$ (see \S~\ref{ss:BB}).
\end{itemize}

The homomorphism 
$$
  \coprod_{\bm{a}} p^{\bm{a}}: \: \check{T}_n = \coprod_{\bm{a}} \check{T}_n^{\bm{a}} \to \Sigma_n
$$
must descend to $p_n: \: T_n \to \Sigma_n$ because $C^*_r(\coprod_{\bm{a}} p^{\bm{a}}): \: C^*_r(G_n, \Sigma_n) \to C^*_r(\check{H}_n, \check{T}_n), \, f \ma f \circ (\coprod_{\bm{a}} p^{\bm{a}})$ lands in $C^*_r(H_n, T_n)$. Moreover, the homomorphisms $\Phi_n$ and $\check{\varphi}_n$ admit similar groupoid models (say $\cP_n$ and $\check{p}_n$) as $\varphi_n$, so that we obtain a commutative diagram
$$
  \xymatrix{
  \cC_n \amalg \cF_n & \ar[l]^(.33){\cP_n} \cT_n \ar[r]_(.3){\subseteq} & \cC_{n+1} \amalg \cF_{n+1}
  \\
  \check{\Sigma}_n \ar[u]_{\subseteq} \ar@{->>}[d] & \ar[l]_{\check{p}_n} \check{T}_n \ar[u]_{\subseteq} \ar[r]^{\subseteq} \ar@{->>}[d] & \check{\Sigma}_{n+1} \ar[u]^{\subseteq} \ar@{->>}[d]
  \\
  \Sigma_n & \ar[l]_{p_n} T_n \ar[r]^{\subseteq} & \Sigma_{n+1}.
  }
$$

\section{Conclusions}

\bproof[Proofs of Theorems~\ref{THM1} and \ref{THM2}]

All we have to do is to check the conditions in Theorem~\ref{THM:CartanLim}, using Proposition~\ref{P:phi=ip} and the groupoid models in \S~\ref{s:GPDModels}. We treat the unital and stably projectionless cases simultaneously. Given a prescribed Elliott invariant, let $A_n$ and $\varphi_n$ be as in \S~\ref{ss:SummConst}. Consider the groupoid models for $A_n$ and $\varphi_n$ in \S~\ref{s:GPDModels}. First of all, by construction, $(H_n,T_n)$ is a subgroupoid of $(G_{n+1},\Sigma_{n+1})$ and $H_n \subseteq G_{n+1}$ is open. Let $(i_n, \imath_n)$ be the canonical inclusion. Secondly, $p_n$ is proper because all the $p^{\bm{a}}$ in \S~\ref{ss:ConnMaps} are proper (they are closed, and pre-images of points are compact). Thirdly, $p_n$ is fibrewise bijective because this is true for $\check{p}_n$ and the canonical projections $\check{\Sigma}_n \onto \Sigma_n$, $\check{T}_n \onto T_n$. By construction, all the connecting maps $\varphi_n$ in \S~\ref{ss:SummConst} are of the form $\varphi_n = (\imath_n)_* \circ (p_n)^*$. Thus, by Proposition~\ref{P:phi=ip}, the conditions in Theorem~\ref{THM:CartanLim} are satisfied. Hence $\ilim \gekl{B_n; \varphi_n}$, with $B_n$ as in \S~\ref{ss:BB}, is a Cartan subalgebra of $\ilim \gekl{A_n; \varphi_n}$, and actually even a C*-diagonal by Remark~\ref{R:Weyl,full,diag} because all $G_n$ are principal.
\eproof

\bremark
\label{R:G0=invlim}
By Remark~\ref{R:Weyl,full,diag}, the twisted groupoids $(G,\Sigma)$ we obtain in the proofs of Theorems~\ref{THM1} and \ref{THM2} are given by the Weyl twists described by \eqref{e:ilimGSigma} and \eqref{e:barSigma,G}. Moreover, it is easy to see that for the groupoids in \S~\ref{s:GPDModels}, we have $(\cC_n \amalg \cF_n)^{(0)} = \check{G}_n^{(0)}$, and since $\Phi_n$, $\check{\varphi}_n$ and $\varphi_n$ send full elements to full elements, $\dot{\cT}_n^{(0)} = (\cC_n \amalg \cF_n)^{(0)}$, $\check{H}_n^{(0)} = \check{G}_{n+1}^{(0)}$ and $H_n^{(0)} = G_{n+1}^{(0)}$, for all $n$ (by Remark~\ref{R:full}). So Remark~\ref{R:Weyl,full,diag} tells us that $G^{(0)} \cong \plim \{ G_n^{(0)}; \dot{p}_n \}$.
\eremark

We now turn to the additional statements in \S~\ref{s:Intro}. In order to prove Corollaries~\ref{C:unbdd1} and \ref{C:unbdd2}, we need to show the following statement. In both the unital and stably projectionless cases, let $A = C^*_r(G,\Sigma)$, $D = C_0(G^{(0)})$, and $\gamma = \ti{\gamma} \vert_T$ be as in Corollaries~\ref{C:unbdd1} and \ref{C:unbdd2}. Let $\cC$ be the canonical diagonal subalgebra of the algebra of compact operators $\cK$.

\bprop
\label{P:da}
There exists a positive element $a \in D \otimes \cC \subseteq A \otimes \cK$ such that $d_{\bullet}(a) = \gamma$.
\eprop
\setlength{\parindent}{0cm} \setlength{\parskip}{0cm}

Here $d_{\bullet}(a)$ denotes the function $T \ni \tau \ma d_{\tau}(a)$. For the proof, we need the following
\blemma
\label{L:dz}
Given a continuous affine map $g: \: T \to (0,\infty)$ and $\varepsilon > 0$, there exists $z \in D \otimes D_k \subseteq A \otimes M_k \subseteq A \otimes \cK$ with $\norm{z} = 1$, $z \geq 0$, $z \in {\rm Ped}(A \otimes \cK)$ such that $g - \varepsilon < d_{\bullet}(z) < g + \varepsilon$.
\elemma
Here $D_k$ is the canonical diagonal subalgebra of $M_k$.
\bproof
We treat the unital and stably projectionless cases simultaneously. Let $\hat{F}_n$ be as in \S~\ref{ss:UnitalCase} and $\hat{F} \defeq \ilim \hat{F}_n$. Choose $a \in \hat{F} \otimes \cK$ with $a \geq 0$ and $d_{\bullet}(a) = g$. Then we can choose $b \in \hat{F}_n \otimes M_k$ (for $n$ big enough) with $b \geq 0$, $d_{\bullet}(b)$ continuous and
$$
  g - \varepsilon < d_{\bullet}(b) < g + \varepsilon.
$$
Using \cite[Theorem~3.1]{APS} just as in \cite[Proof of (6.2) and (6.3)]{Rob}, choose $c \in D(\hat{F}_n) \otimes D_k$ with $c \geq 0$ such that $c$ and $b$ are Cuntz equivalent, where $D(\hat{F}_n)$ is the canonical diagonal subalgebra of $\hat{F}_n$. Choose $d \in {\rm Ped}(A_n \otimes M_k)$ with $d \in D(A_n) \otimes D_k$ such that $(\pi \otimes \id)(d) = c $, where $\pi: \: A_n \onto F_n \onto \hat{F}_n$ is the canonical projection. Let $z$ denote the image of $d$ under the canonical map $A_n \otimes M_k \to A \otimes M_k$. Then $z \in D \otimes D_k$. It is now straightforward to check, using the isomorphism $T(A) \cong T(\hat{F})$ from \cite[\S~13]{GLN}, that $z$ has the desired properties.
\eproof

\bproof[Proof of Proposition~\ref{P:da}]
There is a sequence $(\gamma_i)$ of continuous affine maps $T \to [0,\infty)$ with $\gamma_i \nearrow \gamma - \min(\gamma)$. Choose $\varepsilon_i > 0$ such that $\sum_i \varepsilon_i = \min(\gamma)$. Define $f_i \defeq \gamma_i + \sum_{h=1}^{i-1} \varepsilon_h$. Then $f_i \nearrow \gamma$ and $f_i > 0$. Moreover,
$$
  f_{i+1} = \gamma_{i+1} + \sum_{h=1}^i \varepsilon_h \geq \gamma_i + \Big( \sum_{h=1}^{i-1} \varepsilon_h \Big) + \varepsilon_i = f_i + \varepsilon_i.
$$
Using Lemma~\ref{L:dz}, proceed inductively on $i$ to find $z_i \in D \otimes D_{k(i)}$ such that
$$
  \Big( f_{i+1} - \sum_{h=1}^i d_{\bullet}(z_h) \Big) - \varepsilon_{i+1} < d_{\bullet}(z_{i+1}) < \Big( f_{i+1} - \sum_{h=1}^i d_{\bullet}(z_h) \Big) + \varepsilon_{i+1}.
$$
Note that $f_{i+1} - \sum_{h=1}^i d_{\bullet}(z_h) > 0$ since $\sum_{h=1}^i d_{\bullet}(z_h) < f_i + \varepsilon_i \leq f_{i+1}$.
\setlength{\parindent}{0cm} \setlength{\parskip}{0.5cm}

By construction, we have
$$
  f_i - \varepsilon_i < \sum_{h=1}^i d_{\bullet}(z_h) < f_i + \varepsilon_i, 
  \ 
  \text{so that} \ \sum_{h=1}^i d_{\bullet}(z_h) \nearrow \gamma.
$$
Now set
$$
  a \defeq \sum_{h=1}^{\infty}{}^{\oplus} \ \ 2^{-h} z_h, 
  \
  \text{where we put the elements} \ 2^{-h} z_h \ \text{on the diagonal in} \ D \otimes \cC. 
$$
In this way, we obtain an element $a \in D \otimes \cC \subseteq A \otimes \cK$ with $d_{\bullet}(a) = \gamma$.
\eproof

\bproof[Proof of Corollaries~\ref{C:unbdd1} and \ref{C:unbdd2}]
Given $\ti{\gamma}$ as in Corollaries~\ref{C:unbdd1} and \ref{C:unbdd2}, let $\gamma = \ti{\gamma} \vert_T$. Using Proposition~\ref{P:da}, choose a positive element $a \in D \otimes \cC$ with $d_{\bullet}(a) = \gamma$. In the unital case, it is straightforward to check that we can always arrange $a$ to be purely positive. Then it is straightforward to check that $(\overline{a (A \otimes \cK) a}, \overline{a (D \otimes \cC) a})$ is a Cartan pair. Hence, by \cite[Theorem~5.9]{Ren}, there is a twisted groupoid $(\ti{G},\ti{\Sigma})$ such that $(C^*_r(\ti{G},\ti{\Sigma}), C_0(\ti{G}^{(0)})) \cong (\overline{a (A \otimes \cK) a}, \overline{a (D \otimes \cC) a})$. It is now easy to see (compare also \cite[Corollary~6.12]{GLII}) that $(\ti{G},\ti{\Sigma})$ has all the desired properties.
\eproof

\bproof[Proofs of Corollaries~\ref{C:dim} and \ref{C:ZWZ}]
(i) follows from the observation that we only need the twist if $G_0$ has torsion. The claims in (ii)--(iv) about subhomogeneous building blocks and their spectra follow immediately from our constructions (see also Remark~\ref{R:S3S1}). Moreover, the inverse limit description of the unit space in Remark~\ref{R:G0=invlim} and the dimension formula for inverse limits (see for instance \cite[Chapter~3, \S~5.3, Theorem~22]{Fed}) imply that $\dim(G^{(0)}) \leq 3$ in (ii), $\dim(G^{(0)}) \leq 2$ in (iii) and $\dim(G^{(0)}) \leq 1$ in (iv) and (v). Since $C_0(G^{(0)})$ is projectionless in Theorem~\ref{THM2}, we obtain $\dim(G^{(0)}) \neq 0$, which forces $\dim(G^{(0)}) = 1$ in (v), in the situation of Theorem~\ref{THM2}. In particular, this shows that $\cW$ and $\cZ_0$ have C*-diagonals with one-dimensional spectra. Similarly, given a groupoid $G$ with $\cZ \cong C^*_r(G)$, the only projections in $C(G^{(0)})$ are $0$ and $1$, so that $\dim(G^{(0)}) \neq 0$ and hence $\dim(G^{(0)}) = 1$. It remains to prove that $\dim(G^{(0)}) \geq 3$ in (ii), $\dim(G^{(0)}) \geq 2$ in (iii) and $\dim(G^{(0)}) \geq 1$ in (iv).
\setlength{\parindent}{0cm} \setlength{\parskip}{0.5cm}

To do so, let us use the same notation as in \S~\ref{s:GPDModels}, and write $X_n \defeq G_n^{(0)}$, $Q_n \defeq \dot{\cC}_n^{(0)}$, and $W_n \defeq \dot{\cF}_n^{(0)}$. Clearly, $Q_n$ is homotopy equivalent to a finite set of points, so that for any cohomology theory $H^{\bullet}$ (satisfying the Eilenberg-Steenrod axioms, see \cite[Chapter~17]{tDi}), we have
\begin{equation}
\label{e:H^iQ}
  H^{\bullet}(Q_n) \cong \gekl{0} \ {\rm whenever} \ \bullet \geq 1. 
\end{equation}
Let $P_n \defeq \menge{(t,x) \in Q_n}{t \in \gekl{0,1}, \, x \in \cR_t^{(0)}}$. Then we have a pushout diagram
$$
  \xymatrix{
  P_n \ar[r] \ar[d] & W_n \ar[d] \\
  Q_n \ar[r] & X_n,
  }
$$
where $P_n \to W_n$ is induced by $\bm{b}_t$ and the left vertical arrow is the canonical inclusion. The long exact (Mayer-Vietoris type) sequence attached to the pushout reads
$$
  \dotso \to H^{\bullet - 1}(P_n) \to H^{\bullet}(X_n) \to H^{\bullet}(Q_n) \times H^{\bullet}(W_n) \to H^{\bullet}(P_n) \to H^{\bullet + 1}(X_n) \to \dotso.
$$
Since $H^{\bullet}(P_n) \cong \gekl{0}$ and $H^{\bullet}(Q_n) \cong \gekl{0}$ (see \eqref{e:H^iQ}), we deduce that the canonical map $W_n \to X_n$ induces a surjection $H^{\bullet}(X_n) \to H^{\bullet}(W_n)$ for $\bullet \geq 1$. Moreover, the map
$$
  Q_{n+1} \amalg W_{n+1} = (\dot{\cC}_{n+1} \amalg \dot{\cF}_{n+1})^{(0)} = \check{G}_{n+1}^{(0)} = \check{H}_n^{(0)} \overset{\hat{p}_n}{\lori} \hat{G}_n^{(0)} = Q_n \amalg W_n
$$
induces for $\bullet \geq 1$ a homomorphism $H^{\bullet}(\check{p}_n): \: H^{\bullet}(W_n) \to H^{\bullet}(W_{n+1})$ which fits into the commutative diagram
$$
  \xymatrix{
  H^{\bullet}(X_n) \ar@{->>}[r] \ar[d]_{H^{\bullet}(p_n)} & H^{\bullet}(W_n) \ar[d]^{H^{\bullet}(\check{p}_n)} \\
  H^{\bullet}(X_{n+1}) \ar@{->>}[r] & H^{\bullet}(W_{n+1}).
  }
$$
Thus the canonical maps $W_n \to X_n$ induce for all $\bullet \geq 1$ surjections
$$
  \check{H}^{\bullet}(G^{(0)}) \cong \ilim \gekl{H^{\bullet}(X_n); H^{\bullet}(p_n)} \onto \ilim \gekl{H^{\bullet}(W_n); H^{\bullet}(\check{p}_n)}.
$$
Here $\check{H}^{\bullet}$ is \v{C}ech cohomology, and the first identification follows from the inverse limit description of $G^{(0)}$ in Remark~\ref{R:G0=invlim} and continuity of \v{C}ech cohomology. By construction, $W_n = Z_n \times I_n$ for some finite set $I_n$ and $Z_n$ is as in \S~\ref{ss:SummConst}. Now it is an immediate consequence of our construction that $\ilim \gekl{H^{\bullet}(W_n); H^{\bullet}(\check{p}_n)}$ surjects onto ${\rm Tor}(G_1)$ in case (ii) for $\bullet = 3$, ${\rm Tor}(G_0)$ in case (iii) for $\bullet = 2$, and $G_1$ in case (iv) for $\bullet = 1$. Hence it follows that $\check{H}^3(G^{(0)}) \ncong \gekl{0}$ in case (ii), $\check{H}^2(G^{(0)}) \ncong \gekl{0}$ in case (iii), and $\check{H}^1(G^{(0)}) \ncong \gekl{0}$ in case (iv). As cohomological dimension is always a lower bound for covering dimension, this implies $\dim(G^{(0)}) \geq 3$ in case (ii), $\dim(G^{(0)}) \geq 2$ in case (iii), and $\dim(G^{(0)}) \geq 1$ in case (iv), as desired.
\eproof
\setlength{\parindent}{0cm} \setlength{\parskip}{0.5cm}

\section{Examples}
\label{s:ZWZ}

Let us describe concrete groupoid models for the Jiang-Su algebra $\cZ$, the Razak-Jacelon algebra $\cW$ and the stably projectionless version $\cZ_0$ of the Jiang-Su algebra as in \cite[Definition~7.1]{GLII}. These C*-algebras can be constructed in a way which fits into the framework of \S~\ref{ss:SummConst}, so that our general machinery in \S~\ref{s:LimCartan} produces groupoid models as in \S~\ref{s:GPDModels}. In the following, we focus on $\cZ$.

First we recall the original construction of $\cZ$ in \cite{JS}. For every $n \in \Nz$, choose natural numbers $p_n$ and $q_n$ such that they are relatively prime, with $p_n \mid p_{n+1}$ and $q_n \mid q_{n+1}$, such that $\frac{p_{n+1}}{p_n} > 2 q_n$ and $\frac{q_{n+1}}{q_n} > 2 p_n$. Then $\cZ = \ilim \gekl{A_n; \varphi_n}$, where $A_n = \menge{(f,a) \in C([0,1], E_n) \oplus F_n}{f(t) = \beta_t(a) \ {\rm for} \ t = 0,1}$, $E_n = M_{p_n} \otimes M_{q_n}$, $F_n = M_{p_n} \oplus M_{q_n}$, $\beta_0: \: M_{p_n} \oplus M_{q_n} \to M_{p_n} \otimes M_{q_n}, \, (x,y) \ma x \otimes 1_{q_n}$, $\beta_1: \: M_{p_n} \oplus M_{q_n} \to M_{p_n} \otimes M_{q_n}, \, (x,y) \ma 1_{p_n} \otimes y$. 

To describe $\varphi_n$ for fixed $n$, let $d_0 \defeq \frac{p_{n+1}}{p_n}$, $d_1 \defeq \frac{q_{n+1}}{q_n}$, $d \defeq d_0 \cdot d_1$, and write $d = l_0 q_{n+1} + r_0$ with $0 \leq r_0 < q_{n+1}$, $d = l_1 p_{n+1} + r_1$ with $0 \leq r_1 < p_{n+1}$. Note that we must have $d_1 \mid r_0$ and $d_0 \mid r_1$. Then 
$$
  \varphi_n(f) = u_{n+1}^* \cdot (f \circ \lambda_y)_{y \in \cY(n)} \cdot u_{n+1}, \ \ \ {\rm where} \ \cY(n) = \gekl{1, \dotsc, d} \ {\rm and} \ 
  \lambda_y(t) =
  \bfa
    \frac{t}{2} & {\rm if} \ 1 \leq y \leq r_0,\\
    \frac{1}{2} & {\rm if} \ r_0 < y \leq d - r_1,\\
    \frac{t+1}{2} & {\rm if} \ d - r_1 < y \leq d.
  \efa
$$
Here we think of $A_n$ as a subalgebra of $C([0,1],E_n)$ via the embedding $A_n \into C([0,1],E_n), \, (f,a) \ma f$. 

To construct groupoid models for building blocks and connecting maps, start with a set $\cX(1)$ with $p_1 \cdot q_1$ elements, and define recursively $\cX(n+1) \defeq \cX(n) \times \cY(n)$. Let $\cR(n)$ be the full equivalence relation on $\cX(n)$. Let $\cR(n,p)$ and $\cR(n,q)$ be the full equivalence relations on finite sets $\cX(n,p)$ and $\cX(n,q)$ with $p_n$ and $q_n$ elements. For $t=0,1$, let $\rho_{n+1,t}$ be the bijections corresponding to conjugation by $u_{n+1}(t)$, which induce $\sigma_{n,t}: \: \cR(n,p) \times \cR(n,q) \cong \cR(n)$ corresponding to conjugation by $v_n(t)$ introduced in Remark~\ref{R:unitaries}. Now set
\begin{eqnarray*}
  \check{G}_n &\defeq& \menge{(t,\gamma) \in [0,1] \times \cR(n)}{\gamma \in \sigma_{n,0}(\cR(n,p) \times \cX(n,q)) \ {\rm if} \ t=0, \, \gamma \in \sigma_{n,1}(\cX(n,p) \times \cR(n,q)) \ {\rm if} \ t=1},\\
  G_n &\defeq& \check{G}_n / { }_{\sim} \ \ \ {\rm where} \ \sim \ {\rm is} \ {\rm given} \ {\rm by} \ (0,\sigma_{n,0}(\gamma,y)) \sim (0,\sigma_{n,0}(\gamma,y')) \ {\rm and} \ (1,\sigma_{n,1}(x,\eta)) \sim (1,\sigma_{n,1}(x',\eta)).
\end{eqnarray*}
Define $\check{p}_n: \: \check{H}_n \to \check{G}_n$ as the restriction of $\cP_n: \: \dot{\cT}_n \defeq [0,1] \times \cR(n) \times \cY(n) \to [0,1] \times \cR(n), \, (t,\gamma,y) \ma (\lambda_y(t),\gamma)$ to $\check{H}_n \defeq \cP_n^{-1}(\check{G}_n)$. Set $H_n \defeq \check{H}_n / { }_{\sim}$ where $\sim$ is the equivalence relation defining $G_{n+1} = \check{G}_{n+1} / { }_{\sim}$. The map $\check{p}_n$ descends to $p_n: \: H_n \to G_n$. The groupoid $G$ with $\cZ \cong C^*_r(G)$ is now given by \eqref{e:ilimGSigma} and \eqref{e:barSigma,G}. As explained in Remark~\ref{R:G0=invlim}, its unit space $X \defeq G^{(0)}$ is given by $X \cong \ilim \gekl{X_n; p_n}$, where $X_n = G_n^{(0)}$.

To further describe $X$, let $\bm{p}_n$ be the set-valued function on $[0,1]$ defined by $\bm{p}_n(s) \defeq \menge{\lambda_y(s)}{y \in \cY(n)}$. We can form the inverse limit
$$
  \bm{X} \defeq \plim \gekl{[0,1]; \bm{p}_n} \defeq \Big\{ (s_n) \in \prod_{n=1}^{\infty} [0,1]: \: s_n \in \bm{p}_n(s_{n+1}) \Big\}.
$$
as in \cite[\S~2.2]{IM}. It is easy to see that $X_n \ma [0,1], \, [(t,x)] \ma t$ gives rise to a continuous surjection $X \onto \bm{X}$ whose fibres are all homeomorphic to the Cantor space. Moreover, $\bm{X}$ is connected and locally path connected. The space $X$ itself is also connected. This follows easily from the construction itself (basically from $\gcd(p_n,q_n) = 1$) and also from abstract reasons because $\cZ$ is unital projectionless. In addition, it is straightforward to check that for particular choices for $\rho_{n,t}$ and hence $\sigma_{n,t}$, our space $X$ becomes locally path connected as well. In that case, it is a one-dimensional Peano continuum.

Every $X_n$ is homotopy equivalent to a finite bouquet of circles. It is then easy to compute K-theory and \v{C}ech (co)homology:
\setlength{\parindent}{0cm} \setlength{\parskip}{0cm}

\begin{eqnarray}
  && K_0(C(X)) = \Zz [1], \ \ \ K_1(C(X)) \cong \bigoplus_{i=1}^{\infty} \Zz;\\
  && \check{H}^{\bullet}(X) \cong
  \bfa
  \Zz & {\rm for} \ \bullet = 0,\\
  \bigoplus_{i=1}^{\infty} \Zz & {\rm for} \ \bullet = 1,\\
  \gekl{0} & {\rm for} \ \bullet \geq 2,\\
  \efa
  \ \ \ {\rm and} \ \ \
  \check{H}_{\bullet}(X) \cong
  \bfa
  \Zz & {\rm for} \ \bullet = 0,\\
  \prod_{i=1}^{\infty} \Zz & {\rm for} \ \bullet = 1,\\
  \gekl{0} & {\rm for} \ \bullet \geq 2.\\
  \efa
\end{eqnarray}
It follows that for choices of $\rho_{n,t}$ and $\sigma_{n,t}$ such that $X$ is locally path connected, $X$ must be shape equivalent to the Hawaiian earring by \cite{DV}. In particular, its first \v{C}ech homotopy group is isomorphic to the one of the Hawaiian earring, which is the canonical projective limit of non-abelian free groups of finite rank. Moreover, by \cite{Eda16}, the singular homology $H_1(X)$ coincides with the singular homology of the Hawaiian earring, which is described in \cite{EK}. We refer the reader to \cite{MS} for more information about shape theory, which is the natural framework to study our space since it is constructed as an inverse limit.
\setlength{\parindent}{0cm} \setlength{\parskip}{0.5cm}

Now we turn to $\cW$. Recall the construction in \cite{Jac}. For every $n \in \Nz$, choose integers $a_n, b_n \geq 1$ with $a_{n+1} = 2 a_n + 1$, $b_{n+1} = a_{n+1} \cdot b_n$. Then $\cW = \ilim \gekl{A_n; \varphi_n}$, where $A_n = \menge{(f,a) \in C([0,1], E_n) \oplus F_n}{f(t) = \beta_t(a) \ {\rm for} \ t = 0,1}$, $E_n = M_{(a_n + 1) \cdot b_n}$, $F_n = M_{b_n}$, with 
$$
  \beta_0: \: M_{b_n} \to M_{(a_n + 1) \cdot b_n}, \,
  x \ma 
      \begin{pmatrix}
    x & & & \\
    & \ddots & & \\
    & & x & \\
    & & & 0
    \end{pmatrix}
    \ \ \ {\rm and} \ \ 
  \beta_1: \: M_{b_n} \to M_{(a_n + 1) \cdot b_n}, \,
  x \ma 
      \begin{pmatrix}
    x & & \\
    & \ddots & \\
    & & x 
      \end{pmatrix},
$$
where we put $a_n$ copies of $x$ on the diagonal for $\beta_0$, and $a_n + 1$ copies of $x$ on the diagonal for $\beta_1$. To describe $\varphi_n$ for fixed $n$, let $d \defeq 2 a_{n+1}$. Then 
$$
  \varphi_n(f) = u_{n+1}^* \cdot (f \circ \lambda_y)_{y \in \cY(n)} \cdot u_{n+1}, \ \ \ {\rm where} \ \cY(n) = \gekl{1, \dotsc, d} \ {\rm and} \ 
  \lambda_y(t) =
  \bfa
    \frac{t}{2} & {\rm if} \ 1 \leq y \leq a_{n+1},\\
    \frac{1}{2} & {\rm if} \ y = a_{n+1} + 1,\\
    \frac{t+1}{2} & {\rm if} \ a_{n+1} + 1 < y \leq d.
  \efa
$$
Here we think of $A_n$ as a subalgebra of $C([0,1],E_n)$ via the embedding $A_n \into C([0,1],E_n), \, (f,a) \ma f$.

To construct groupoid models, start with a set $\cX(1)$ with $(a_1 + 1) \cdot b_1$ elements, and define recursively $\cX(n+1) \defeq \cX(n) \times \cY(n)$. Let $\cR(n)$ be the full equivalence relation on $\cX(n)$. Let $\cR(n,a)$ and $\cR(n,b)$ be the full equivalence relations on finite sets $\cX(n,a)$ and $\cX(n,b)$ with $a_n + 1$ and $b_n$ elements, and let $\cX'(n,a) \subseteq \cX(n,a)$ be a subset with $a_n$ elements (corresponding to the multiplicity of $\beta_0$). For $t=0,1$, let $\rho_{n+1,t}$ be the bijections corresponding to conjugation by $u_{n+1}(t)$, which induce $\sigma_{n,t}: \: \cR(n,a) \times \cR(n,b) \cong \cR(n)$ corresponding to conjugation by $v_n(t)$ introduced in Remark~\ref{R:unitaries}. Set
\begin{eqnarray*}
  \check{G}_n &\defeq& \menge{(t,\gamma) \in [0,1] \times \cR(n)}{\gamma \in \sigma_{n,0}(\cX'(n,a) \times \cR(n,b)) \ {\rm if} \ t=0, \, \gamma \in \sigma_{n,1}(\cX(n,a) \times \cR(n,b)) \ {\rm if} \ t=1},\\
  G_n &\defeq& \check{G}_n / { }_{\sim} \ \ \ {\rm where} \ \sim \ {\rm is} \ {\rm given} \ {\rm by} \ (t,\sigma_{n,t}(x,\gamma)) \sim (t',\sigma_{n,t'}(x',\gamma)).
\end{eqnarray*}
Now define $\check{p}_n: \: \check{H}_n \to \check{G}_n$ as the restriction of $\cP_n: \: \dot{\cT}_n \defeq [0,1] \times \cR(n) \times \cY(n) \to [0,1] \times \cR(n), \, (t,\gamma,y) \ma (\lambda_y(t),\gamma)$ to $\check{H}_n \defeq \cP_n^{-1}(\check{G}_n)$. Set $H_n \defeq \check{H}_n / { }_{\sim}$ where $\sim$ is the equivalence relation defining $G_{n+1} = \check{G}_{n+1} / { }_{\sim}$. The map $\check{p}_n$ descends to $p_n: \: H_n \to G_n$. The groupoid $G$ with $\cW \cong C^*_r(G)$ is now given by \eqref{e:ilimGSigma} and \eqref{e:barSigma,G}. As explained in Remark~\ref{R:G0=invlim}, its unit space $X \defeq G^{(0)}$ is given by $X \cong \plim \gekl{X_n; p_n}$, where $X_n = G_n^{(0)}$. As in the case of $\cZ$, $X$ surjects continuously onto $\plim \gekl{\Tz; \bm{p}_n}$ with Cantor space fibres, where $\Tz = [0,1] / { }_{0 \sim 1}$ and $\bm{p}_n([s]) = \menge{[\lambda_y(s)]}{y \in \cY(n)}$. However, it is easy to see that (at least for some choices of $\rho_{n,t}$ and $\sigma_{n,t}$), $X$ will not be connected, though its connected components all have to be non-compact.

Now let us treat $\cZ_0$. For each $m \in \Nz$, choose integers $a_n, b_n, h_n \geq 1$ with $a_{n+1} = ((2a_n + 2)h_n + 1) \cdot a_n$, $b_{n+1} = ((2a_n + 2)h_n + 1) \cdot b_n$. Let $A_n = \menge{(f,a) \in C([0,1], E_n) \oplus F_n}{f(t) = \beta_t(a) \ {\rm for} \ t = 0,1}$, with $E_n = M_{(2a_n + 2) \cdot b_n}$, $F_n = M_{b_n} \oplus M_{b_n}$,
$$
  \beta_0: \: F_n \to E_n, \,
  (x,y) \ma
  \rukl{
    \begin{smallmatrix}
    x & & & & & & &\\
    & \ddots & & & & & &\\
    & & x & & & & &\\
    & & & 0 & & & &\\
    & & & & y & & &\\
    & & & & & \ddots & &\\
    & & & & & & y &\\
    & & & & & & & 0
    \end{smallmatrix}
  },
  \ \ \ {\rm and} \ \ \
  \beta_1: \: F_n \to E_n, \,
  (x,y) \ma 
  \rukl{
    \begin{smallmatrix}
    x & & & & &\\
    & \ddots & & & &\\
    & & x & & &\\
    & & & y & &\\
    & & & & \ddots &\\
    & & & & & y\\
    \end{smallmatrix}
  },
$$
where we put $a_n$ copies of $x$ and $y$ on the diagonal for $\beta_0$, and $a_n + 1$ copies of $x$ and $y$ on the diagonal for $\beta_1$. To describe the connecting maps $\varphi_n: \: A_n \to A_{n+1}$, fix $n$ and let $d \defeq (2 a_{n+1} + 2)h_n + (2 a_n h_n + 1)$. Then $(2a_{n+1} + 2) \cdot b_{n+1} = d \cdot (2a_n + 2) \cdot b_n$. It is now easy to see that for suitable choices of unitaries $u_{n+1}$, whose values at $0$ and $1$ are permutation matrices, we obtain a homomorphism $\varphi_n: \: A_n \to A_{n+1}$ by setting
$$
  \varphi_n(f) \defeq u_{n+1}^* \cdot (f \circ \lambda_y)_{y \in \cY(n)} \cdot u_{n+1}, \ {\rm for} \ \cY(n) = \gekl{1, \dotsc, d}, \, 
  \lambda_y(t) =
  \bfa
    \frac{t}{2} & {\rm if} \ 1 \leq y \leq 2 a_k h_k + 2h_k + 1,\\
    \frac{1}{2} & {\rm if} \ 2 a_k h_k + 2h_k + 1 < y \leq (2a_{k+1} + 2) h_k,\\
    \frac{t+1}{2} & {\rm if} \ (2a_{k+1} + 2) h_k < y \leq d.
  \efa
$$
As above, we think of $A_n$ as a subalgebra of $C([0,1],E_n)$ via $A_n \into C([0,1],E_n), \, (f,a) \ma f$. Now arguments similar to those in \cite{JS,Jac} show that $\ilim \gekl{A_n; \varphi_n}$ has the same Elliott invariant as $\cZ_0$, so that $\cZ_0 \cong \ilim \gekl{A_n; \varphi_n}$ by \cite[Corollary~6.2.4]{Rob} (see also \cite[Theorem~12.2]{GLII}).

To construct groupoid models, start with a set $\cX(1)$ with $(2a_1 + 2) \cdot b_1$ elements, and define recursively $\cX(n+1) \defeq \cX(n) \times \cY(n)$. Let $\cR(n)$ be the full equivalence relation on $\cX(n)$. Let $\cR(n,a,1)$, $\cR(n,a,2)$, $\cR(n,b,1)$ and $\cR(n,b,2)$ be full equivalence relations on finite sets $\cX(n,a,1)$, $\cX(n,a,2)$, $\cX(n,b,1)$ and $\cX(n,b,2)$ with $a_n + 1$, $a_n + 1$, $b_n$ and $b_n$ elements, respectively. Let $\cX_0(n,a,1) \subseteq \cX(n,a,1)$ and $\cX_0(n,a,2) \subseteq \cX(n,a,2)$ be subsets with $a_n$ elements (corresponding to the multiplicities of $\beta_0$), and set $\cX_1(n,a,\bullet) \defeq \cX(n,a,\bullet)$. For $t=0,1$, let $\rho_{n+1,t}$ be the bijections corresponding to conjugation by $u_{n+1}(t)$, which induce $\sigma_{n,t}: \: \cR(n,a,1) \times \cR(n,b,1) \amalg \cR(n,a,2) \times \cR(n,b,2) \cong \cR(n)$ corresponding to conjugation by $v_n(t)$ introduced in Remark~\ref{R:unitaries}. Set
\begin{eqnarray*}
  \check{G}_n &\defeq& \menge{(t,\gamma) \in [0,1] \times \cR(n)}{\gamma \in \sigma_{n,t}(\cX_t(n,a,1) \times \cR(n,b,1) \amalg \cX_t(n,a,2) \times \cR(n,b,2)) \ {\rm if} \ t \in \gekl{0,1}},\\
  G_n &\defeq& \check{G}_n / { }_{\sim} \ \ \ {\rm where} \ \sim \ {\rm is} \ {\rm given} \ {\rm by} \ (t,\sigma_{n,t}(x,\gamma)) \sim (t',\sigma_{n,t'}(x',\gamma)).
\end{eqnarray*}
Now define $\check{p}_n: \: \check{H}_n \to \check{G}_n$ as the restriction of $\cP_n: \: \dot{\cT}_n \defeq [0,1] \times \cR(n) \times \cY(n) \to [0,1] \times \cR(n), \, (t,\gamma,y) \ma (\lambda_y(t),\gamma)$ to $\check{H}_n \defeq \cP_n^{-1}(\check{G}_n)$. Set $H_n \defeq \check{H}_n / { }_{\sim}$ where $\sim$ is the equivalence relation defining $G_{n+1} = \check{G}_{n+1} / { }_{\sim}$. The map $\check{p}_n$ descends to $p_n: \: H_n \to G_n$. The groupoid $G$ with $\cZ_0 \cong C^*_r(G)$ is now given by \eqref{e:ilimGSigma} and \eqref{e:barSigma,G}. As explained in Remark~\ref{R:G0=invlim}, its unit space $X \defeq G^{(0)}$ is given by $X \cong \plim \gekl{X_n; p_n}$, where $X_n = G_n^{(0)}$. As for $\cW$, $X$ surjects continuously onto $\plim \gekl{\Tz; \bm{p}_n}$ with Cantor space fibres, where $\Tz = [0,1] / { }_{0 \sim 1}$ and $\bm{p}_n([s]) = \menge{[\lambda_y(s)]}{y \in \cY(n)}$. However, it is easy to see that (at least for some choices of $\rho_{n,t}$ and $\sigma_{n,t}$), $X$ will not be connected, though its connected components all have to be non-compact.

\end{document}